\documentclass[twoside,12pt]{article}

\usepackage[margin=.8in]{geometry}

\usepackage[ruled, vlined, linesnumbered]{algorithm2e}
\usepackage{algorithmic}
\usepackage{subcaption}
\usepackage{color}
\usepackage{bbm,bm}

\usepackage{crop}%
\usepackage{amsmath,amssymb,amsfonts,upref,endnotes}
\usepackage{amsthm}
\usepackage{marginnote}
\usepackage{soul}
\RequirePackage{graphicx}
\usepackage[figuresright]{rotating}
\usepackage{color}

\usepackage{url}
\usepackage{hyperref}
\hypersetup{
	colorlinks=true,
	linkcolor=blue,
	filecolor=magenta,      
	urlcolor=cyan,
}

\usepackage[]{footmisc}

\usepackage{authblk}

\newcommand{\nonl}{\renewcommand{\nl}{\let\nl\oldnl}}

\def\*#1{\mathbf{#1}}
\def\+#1{\mathcal{#1}}
\def\R{\mathbb{R}}

\def\Expect{\mathbb{E}}

\def\vecs{\mathrm{vec}}

\renewcommand{\tilde}{\widetilde}

\newtheorem{theorem}{\bf Theorem}[section]

\newtheorem{corollary}{\bf Corollary}[section]
\newtheorem{definition}{\bf Definition}[section]

\newtheorem{assumption}{\bf Assumption}[section]

\begin{document}

\title{Regularization by Denoising Sub-sampled Newton Method for Spectral CT Multi-Material Decomposition}

\author{Alessandro Perelli}
\author{Martin S. Andersen}

\affil{Technical University of Denmark (DTU Compute)}
\date{\vspace{-1cm}}                     
\renewcommand\Affilfont{\small}

\maketitle

\begin{abstract}
Spectral Computed Tomography (CT) is an emerging technology that enables to estimate the concentration of basis materials within a scanned object by exploiting different photon energy spectra. In this work, we aim at efficiently solving a model-based maximum-a-posterior problem to reconstruct multi-materials images with application to spectral CT. In particular, we propose to solve a regularized optimization problem based on a plug-in image-denoising function using a randomized second order method. By approximating the Newton step using a sketching of the Hessian of the likelihood function, it is possible to reduce the complexity while retaining the complex prior structure given by the data-driven regularizer. We exploit a non-uniform block sub-sampling of the Hessian with inexact but efficient Conjugate gradient updates that require only Jacobian-vector products for denoising term. Finally, we show numerical and experimental results for spectral CT materials decomposition. 
\end{abstract}

	
\section{Introduction}
Spectral Computed Tomography (CT) is a developing technique for diagnostic and therapy in medical imaging which can gather enhanced physical parameters information compared to conventional CT by exploiting the dependence of the attenuation coefficients of an object with distinct photon energy spectra. 
This information leads to the estimation of images of basis materials \cite{Alvarez1976}, which can be used for quantitative imaging \cite{Alessio2013}, reducing beam hardening artifacts and improving the contrast-to-noise ratio \cite{Yu2011}. Recently, spectral CT with photon counting detectors has attracted increasing interest. 


\maketitle

In contrast to conventional transmission CT, photon counting detectors can resolve the X-ray photon energy and allow to acquire registered energy selective images \cite{Roessl2007} with effective reduction of the electronic readout noise, which is of crucial importance for low-dose applications. However, material decomposition algorithms typically lead
to a degradation of the signal-to-noise ratio (SNR) and noise amplification \cite{kapper2010} compared to the unprocessed spectral images. This is still a fundamental problem of spectral CT since it limits the usability of material selective images. 

Most spectral CT algorithms separate the process of material decomposition and image reconstruction. Material decomposition is either performed in the projection space prior to image reconstruction \cite{Alvarez2011,Vilches2017} or an image based decomposition algorithm is applied after reconstructing the images corresponding to different photon energy spectra \cite{Maab2009,Gao2011,Kim2015,Kazantsev2018} by using tensor factorization \cite{wu2020} or dictionary learning \cite{wu2020image}. Separating these steps is sub-optimal because the full statistical information contained in the spectral tomographic measurements cannot be exploited. Statistical iterative reconstruction (SIR) techniques provide a data consistent approach to obtaining material selective images using a statistical model of the projection measurements and incorporating prior knowledge about the reconstructed material images. Unfortunately, the loss function induced by the forward model which directly maps the expected spectral projection measurements and the material images is non-convex and difficult to minimize; convex approximations of the loss function such as convex-concave techniques \cite{Long2014,Barber2016,Schmidt2017} or semi-empirical models \cite{Weidinger2016, Mechlem2018} have been proposed but suffer of high computational cost because of the double-loop nature of the solver. 

The main challenges which are attracting considerable interest are related to the fact that the image to be recovered lies in an high-dimensional space and to leverage more sophisticated or learned data-dependent priors which are needed since the estimation problem is generally ill-posed. Regarding the issue of reducing the per-iteration computational cost to solve the reconstruction in the high-dimensional setting, numerous first order solvers which leverage random numerical linear algebra have been deployed in different machine learning and imaging applications. Well established algorithms like Stochastic Gradient Descent (SGD) \cite{needell2014, konevcny2015},  order subsets \cite{kim2014} and Alternating Direction Method of Multiplier (ADMM) have gained considerably attention in spectral CT \cite{Sawatzky2014, ramani2011}. Unfortunately, the iteration complexity of first order methods relies significantly on the condition number of the problem like SGD which has a sub-linear convergence rate to only a neighborhood of the solution \cite{gower2019}. 
Instead, the information related to the curvature is decisive to improve the convergence by reducing the condition number and to obtain physical meaningful and accurate quantitative estimation. Therefore, accounting for a greedy selection of parameters and step-size is crucial to design first-order algorithms for imaging like CT image reconstruction \cite{kim2014, he2018}. 
Deterministic second order methods for solving regularized optimization problems with Generalized Linear Models (GLM) have been widely studied \cite{schmidt2011, lee2014} but despite the superior, linear or super-linear, convergence rate of Newton methods compared to first order methods one weakness relies on the high computational cost of calculating the Hessian matrix in high dimensional imaging application \cite{schmidt2011, lee2014}. The Newton method requires at each iteration to solving exactly a linear system with complexity $\+O(mn^2)$, with $m$ the dimension of the measurement vector and $n$ the dimension of the input vector, and therefore it is practical only when the dimension $n$ is low or if the Hessian can be factorized in time linear to the input dimension, like by the fast Fourier Transform (FFT) \cite{nocedal2006}. 
Instead inexact methods solve approximately the linear system by iterative algorithms like Conjugate Gradient (CG) which does not require to calculate the full Hessian matrix but only operate through vectors products \cite{gene1996}; this is particularly appealing in imaging systems like CT where there are fast implementations of the operator that represents the physical forward model.

In this work, a new second-order algorithm for statistical iterative material decomposition (Denoising-IHS) is developed which enjoys an improved accuracy and computation trade-off. We propose a computational efficient quantitative estimation with second order and parameters-free methods based on randomized linear algebra and regularization induced by denoising. We exploit the recent development of randomized sketches \cite{woodruff2014} and sub-sampling methods \cite{bollapragada2019} for the Hessian matrix to improve upon the computational cost-per-iteration. Unfortunately, these methods are efficient in the regime where the measurement dimension $m$ is much larger than the input dimension $n$, i.e. $m\gg n$ \cite{berahas2020}, since the theoretical bound on the embedding of the sketched Hessian matrix requires a number of samples of order $\+O(n)$.
Furthermore, these methods have been mainly analyzed for minimizing smooth and convex cost function and the interesting case of using denoiser-based regularizers has not been deeply studied. Therefore, we aim at developing an algorithm for spectral CT model where the dimension of the measurements is of the same order as the dimension of the image to be reconstructed, i.e. $m\approx n$ and exploit denoising function to capture the complex structure of the material images to be reconstructed.


\subsection{Related works}
 
Important contributions related to sketching methods for solving optimization problems refer to the Newton-Sketch algorithm \cite{pilanci2017} where the Newton equations are solved using an unbiased estimator of the true Hessian obtained by performing a dimensionality reduction using random sketching matrices. Different types of sketching such as Sub-Gaussian matrices, Hadamard or Fourier transforms, sparse Count Sketch matrices can be employed and have been analyzed \cite{pilanci2017}. Furthermore, sub-sampling Newton methods \cite{byrd2011, erdogdu2015, roosta2016, xu2016} are based on random row sub-sampling according to a uniform or non-uniform probability distribution \cite{roosta2016, roosta2016ii,  bollapragada2019}. It is interesting to note that random matrix sub-sampling schemes can be considered a particular type of sketching and they are more efficient especially in imaging problems where the linear forward model has a precise, deterministic structure, like the Random transform for CT imaging. Furthermore, since the Radon matrix $\*A$ is an high dimensional tall and sparse matrix, the count sketch might be efficient because of $\+O(\mathrm{nnz}(\*A))$ multiplication complexity, where $\mathrm{nnz}(\*A)$ stands for the number of non-zero elements of the matrix $\*A$. Unfortunately, the count sketch requires a bound on the sketch size of order $\+O(n)$, with $n$ the dimension of the vectorized image input, which is practically useless in CT scenarios where $n$ is of the same order of the number of measurements $m$. 

Our main interest is to develop a non-uniform sub-sampling Newton algorithm for CT imaging by incorporating data-dependent regularizer which is becoming state-of-the art in model-based iterative reconstruction. In imaging applications where the input prior is difficult to model, powerful regularization techniques are based on data-driven models or by exploiting denoising functions, like the Plug-and-play (PnP) approach \cite{venkatakrishnan2013} or the regularization by denoising (RED) \cite{romano2017, reehorst2019, sun2019}. It is interesting how both approaches allow to integrate in the iterative solver either deterministic denoisers, like BM3D \cite{dabov2007}, NLM, TV or deep learning CNN-based denoisers \cite{zhang2017}. The RED approach is based on an explicit regularization term \cite{romano2017} and scalable first-order solvers for the RED problem in high dimensions have been proposed \cite{sun2019block}. As far as the authors are aware, the RED approach has been investigated with first order solvers but not with Newton methods mainly because by default it would require to explicitly generate the Jacobian matrix of the denoiser which is infeasible in high-dimensions.


\subsection{Main Contribution and Relation to Existing Works}

In this work, we develop a second order method (Denoising-IHS) which combine the iterative sketching \cite{pilanci2017} for the log-likelihood function and the RED regularization which can be implemented by a generic denoiser through the score matching formulation with application to spectral CT material decomposition. We propose an effective choice for both the sketching scheme and the solver in order to reduce the computational complexity for imaging problems characterized by physical-based convolutional operators. 

The overall algorithmic time complexity is $T = n_{outer}\cdot(t_{samp} + t_{solve})$ where $t_{samp}$ is the time needed at each iteration to construct the sketching, $t_{solve}$ is the time to solve the linear sub-problem inexactly, and $n_{outer}$ is the total number of iterations which is affected by the convergence rate. 
Considering $t_{samp}$, taking inspiration from \cite{xu2016}, we develop a random non-uniform sub-sampling of the measurement where the probabilities are calculated according to a penalized block-leverage score metric. While penalizing the ridge leverage scores allows to reduce the sketch size with provable convergence, this metric exploits the CT measurements block structure which is defined as a product of the number of views and the number of detectors. This means that we randomly sub-sample over the acquired  angles and the ridge leverage score associated with each view is calculated over a block constituted by the number of detectors. We take advantage of the convolutional structure of the Radon transform to estimate the leverage scores through the Fourier decomposition and reduce the complexity to $\+O(n\log n)$. Regarding the computational time to solve the sub-problem $t_{solve}$, we employ the CG algorithm since it does not require to store any gradient or Jacobian matrices and all operations can be performed using implicitly operators instead of matrix form. We develop a method to estimate the denoiser's Jacobian product through automatic differentiation for CNN-based denoisers. Finally, for the term $n_{outer}$, we provide convergence analysis of the proposed method.


\subsection{Notation}
Matrices or discrete operators and column vectors are written respectively in capital and normal boldface type, i.e. $\*A$ and $\*a$ to distinguish from scalars and real-valued variables written in normal weight. $(\cdot^T)$ refers to the transpose of a matrix and $\*1$ refers to a vector of ones. Non-random quantities and random realizations are not distinguished typographically while random variable are written with capital letters. The conditional probability density function of $\*y$ given $\*x$ is denoted by $p(\*y|\*x)$. A Gaussian random variable $\*x$ with vector mean $\*a$ and isotropic variance $b$ is denoted by $\*x\sim\mathcal{N}(\*a,b\*I)$. $\langle\*a,\*b\rangle=\*b^T\*a$ refers to the vector inner product.

\subsection{Structure of the Paper}
The remainder of this paper is structured as follows: in Section \ref{sec:Spec_CT}, the spectral CT mathematical model with the Gaussian approximation of the noise is described and Section \ref{sec:opt_prob} analyses the optimization problem together with a review of the Hessian sketching Newton method and the formulation of the regularization by denoising approach. The proposed algorithms is presented in Section \ref{sec:proposed_alg} which contains the details of the denoising-based regularizer, the spectral CT data loss and an analysis of the CG inner solver. The further sections describe the crucial blocks of the Denoising-IHS algorithm: the approach to estimating the denoiser directional derivatives needed in the CG solver is described in Section \ref{sec:dir_der} while the sketching strategy based on block leverage scores sub-sampling of the Hessian matrix of the data loss function is presented in Section \ref{sec:blk_scores}. Finally, in Sections \ref{sec:num_results} and \ref{sec:exp_res} we show the numerical simulation and real acquisition on different sets of images and we assess the algorithm in terms of accuracy and computation.


\section{Spectral X-ray CT Model}\label{sec:Spec_CT}

In this Section, we describe the CT physical measurement process in the continuous domain with the spectrum of the X-ray source beams composed by different energy intervals. The spectral X-ray mathematical model is based on the Beer's law which provides the X-ray intensity after transmission, $p_i^k$ of the $i$-th ray at the $k$-th energy interval (bin) as follows
\begin{equation}\label{eq:spCT_cc}
p^k_i = \int_E S^{k}(E) e^{ -\int_{l\in L_i} \mu (\vec{\kappa}(l), E) dl } dE + \eta_i^{k}, \quad \left\{
\begin{array}{lcl}
i & = & 1, 2, \ldots , N_d \times N_p, \\
k & = & 1, 2, \ldots , N_b,
\end{array}\right.
\end{equation}
where $p_i^k$ is the X-ray intensity of the $i$-th detector pixel in the $k$-th detector bin, $E$ is the photon flux density, $N_d$ is the number of detector pixels, $N_p$ is the number of projections, $\mu(\vec{\kappa}(l), E)$ denotes the linear attenuation coefficient that is related to the position function $\vec{\kappa}(l)$ and the energy level $E$, $\eta_i^k$ is the Gaussian error term for the $i$-th element in the $k$-th energy bin and $S^k(E)$ represents the photon flux density for the $k$-th detector bin, which is the number of incident photons at the energy $E$ in the $k$-th energy window.  

In the discrete domain, we map the continuous spatially and energy dependent distribution $\mu(\vec{\kappa}(l), E)\in L_2(\R)$ by using a basis material decomposition for representing the object

\begin{equation}\label{eq:mat_decomp}
\mu(\vec{\kappa}(l), E) \approx \sum_{m=1}^{N_m}c_m(E)\sum_{j=1}^{N_v} \phi_j(\vec{\kappa})x_{j,m}
\end{equation}
where $c_m(E)$ indicates the linear attenuation coefficient of the $m$-th material at the energy level $E$, $\phi_j(\vec{\kappa})$ defines the $j$-th among $N_v$ basis functions associated with a discrete sampling on a $\sqrt{N_v} \times \sqrt{N_v}$ Cartesian grid, $N_m$ represents the number of basis materials and $x_{j,m}$ is the weight fraction of the $m$-th material in the $j$-th image voxel/pixel. According to the basis material parameterization in Eq. (\ref{eq:mat_decomp}), the line integral in Eq. (\ref{eq:spCT_cc}) becomes a summation:
\begin{equation}\label{eq:linCT}
\int_{l\in L_i}\mu(\vec{\kappa}(l), E)dl \;\; \approx \;\; \sum_{m=1}^{N_m} c_m(E) \sum_{j=1}^{N_v} x_{j,m} \int_{l\in L_i}\phi_j(\vec{\kappa}(l))dl \;\; = \;\; \sum_{m=1}^{N_m} \sum_{j=1}^{N_v} r_{i,j} x_{j,m} c_{m,e}
\end{equation}
\noindent where $r_{i,j}$ represents the $i,j$ element of the system matrix describing the discrete summation along the $i$-path from source through object at pixel position $j$ onto each detector. Repeating this over all lines defines the full view linear tomographic system matrix $\*R = [r_{ij}]$; by physical design, $\*R$ is a sparse and non-negative matrix with dimensions $N_d  N_p \times N_v$. For a discrete set of energies of dimension $N_e$, the discrete model (\ref{eq:spCT_cc}) can be written as
\begin{equation}\label{eq:vect_spmodel}
p^k_i = \sum_{e=1}^{N_e} s^k e^{\sum_{m=1}^{N_m}\sum_{j=1}^{N_v} r_{i,j} x_{j,m} c_{m,e}} + \eta_i^k
\end{equation}

\noindent The source spectrum is modeled as a vector $\*s\in\mathbb{R}^{N_e\times 1}$ with non-negative entries and the detector response is a positive matrix $\*D\in\mathbb{R}^{N_b\times N_e}$ whose dimension is the product of the number of bins $N_b$ and the number of discrete energies $N_e$. Therefore, the spectrum for all the energy bins $b=1,\ldots, N_b$ is the matrix obtained as Hadamard point-wise product $\*S = \*D\odot\*s\in\mathbb{R}^{N_b\times N_e}$. Both the source spectrum and the detectors response matrix are assumed to be known. 

In matrix notation, $\*X = \left[\*x_1,\ldots, \*x_l, \ldots \*x_{N_m}\right]\in\mathbb{R}^{N_v\times N_m}$ denotes the matrix obtained as columns concatenation of the vectorized image $\*x_l$ of the $l=1,\ldots, N_m$ basis materials. Therefore, $r_{i,j}$, and $c_{m,e}$ are collected in a matrix form and $p^k_i$, $s^k_e$, $\eta^k_i$ are concatenated with respect to their energy windows. The corresponding matrix equation of (\ref{eq:vect_spmodel}) can be represented as 
\begin{equation}\label{eq:meas_spect}
\*P = \exp\left({-\*R\*X\*C^T}\right)\*S + \*N
\end{equation}

\noindent with $\*C$ a matrix of the size $N_e\times N_m$ where each entry corresponds to $c_{e,m}$, the linear attenuation coefficient of the energy $e$, and the $m$-th material. We vectorize the nonlinear matrix equation (\ref{eq:meas_spect}) on both sides and linearize to avoid a tensor form for computing second order derivatives and we follow the mathematical model derived in \cite{Hu2019}. In the forward problem, we use the full spectrum, and the matrix $\*S$ is rectangular while for solving the inverse problem, we choose the average in each energy window to represent the corresponding energy spectrum, so $N_b = N_e$ and the matrix $\*S$ in is a nonsingular diagonal matrix. Therefore, we multiply $\*S^{-1}$ and vectorize on both sides of (\ref{eq:meas_spect}) and using the properties of Kronecker products we obtain 
\begin{equation}\label{eq:pr_meas_spect}
\left( \*S^{-T}\otimes \*I\right)\*p = \exp\Big({-\left(\*C\otimes\*R\right)\*x}\Big) + \left( \*S^{-T}\otimes \*I\right)\*n
\end{equation}

\noindent where $\*p = \vecs\left(\*P\right)$, $\*x = \vecs\left(\*X\right)$ and $\*n = \vecs\left(\*N\right)$. Denoting $\tilde{\*p} = \left( \*S^{-T}\otimes \*I\right)\*p$ and $\tilde{\*n} = \left( \*S^{-T}\otimes \*I\right)\*n$ and taking the logarithm, we have $\log\left(\tilde{\*p} - \tilde{\*n}\right) = -\left(\*C\otimes\*R\right)\*x$. 

\noindent Using a first order Taylor expansion at $\tilde{\*p}$, $\log \left( \tilde{\*p} - \tilde{\*n} \right) = \log\left(\tilde{\*p}\right) - \mathrm{diag}\left(\tilde{\*p}\right)^{-1}\tilde{\*n} + \+O\left(\|\tilde{\*n}\|^2\right)$, we obtain the following approximation $\*y\approx \left(\*C\otimes\*R\right)\*x -  \mathrm{diag}\left(\tilde{\*p}\right)^{-1}\tilde{\*n}$ where $\*y = -\log\left(\tilde{\*p}\right)$. Assuming Gaussian noise $\*n\sim \+N \big(\*0, \mathrm{diag}(\*p)  \big)$, we have
\begin{equation}
\*y|\*x \sim \+N\big(\left(\*C\otimes\*R\right)\*x, \bm\Sigma \big)
\end{equation}

\noindent where the inverse covariance matrix $\bm\Sigma^{-1}$ is expressed as

\begin{equation}\label{eq:inv_var}
\bm\Sigma^{-1} = \mathrm{diag}\left(\tilde{\*p}\right)\Big( \*S\otimes \*I\Big)\mathrm{diag}\left({\*p}\right)^{-1}\left( \*S^T\otimes \*I\right)\mathrm{diag}\left(\tilde{\*p}\right)
\end{equation}

\noindent $\*P$ is a matrix which represents the number of photons of each energy window in the corresponding column and each element of $\*P$ is a positive integer. From expression (\ref{eq:inv_var}), we can see that the $\bm\Sigma^{-1}$ is diagonal since we have assumed that $\*S$ is diagonal. 

We consider the MAP estimation problem for the spectral X-ray CT model with Gaussian approximation. Given the measurement vector $\*y\in\R^m$ of dimension $m=N_d\cdot N_p\cdot N_b$, we consider the data-fit loss function $f(\*x,\*y):\R^n\times \R^m\rightarrow \R$ as a measure of the distance between the input parameter vector $\*x\in\R^n$, of dimension $n=N_v\cdot N_m$, and the measurement vector $\*y$, with $m > n$. The loss function is given by the following empirical negative log-likelihood (NLL)
\begin{equation}\label{eq:NLPoisson}
f(\*x,\*y) \;\; = \;\; -\sum_{i=1}^m \log p(y_i|\*x) \;\; \approx \;\; \frac{1}{2}\sum_{i=1}^m \sigma_i^{-1}\Big( \*a_i\*x - \*y_i \Big)^2 
\;\; = \;\; \frac{1}{2}\|\*A\*x - \*y\|^2_{\bm\Sigma^{-1}}
\end{equation}

\noindent where $\*A = \*C\otimes\*R$. Given the likelihood function in Eq. (\ref{eq:NLPoisson}) with $\*x\in\R_+^n$, it results that
\begin{equation}\label{eq:HessianCT}
\nabla_{\*x} f(\*x,\*y) \;\;= \;\; \*A^T\bm\Sigma^{-1}\left( \*A\*x - \*y \right),\; \quad\quad\; \nabla^2_{\*x} f(\*x,\*y) \;\;= \;\;\*A^T \bm\Sigma^{-1}\*A   
\end{equation}

\noindent It is important to note that in the region of physical interest, i.e. $a_i^T\*x\in\R_+$, $\nabla^2_{\*x} f(\*x,\*y)$ is bounded.


\section{Optimization problem}\label{sec:opt_prob}

The MAP estimator can be formulated as the following optimization problem 
\begin{equation}\label{eq:pr_GLM}
\hat{\*x} = \arg\min_{\*x\in\R^n_+} f(\*x,\*y) + \rho_{\bm\sigma, \nu}(\*x) 
\end{equation}

\noindent where $f(\*x,\*y) = -\log p(\*y|\*x)$ represents the likelihood loss function as described in Eq. (\ref{eq:NLPoisson}) and $\rho_{\bm\sigma, \nu}(\*x)$ is the regularization term which defines a prior $\rho_{\bm\sigma, \nu}(\*x) = -\log p(\*x, \nu)$ on the vector $\*x$ to be estimated, with $\nu$ the noise level. In the following, we will consider the constrained minimization over the set of non-negative vectors $\*x\in\R^n_+$. We focus on a particular non-convex regularizer which can be explicitly defined according to a denoising (RED) function dependent on a set of parameters $\bm\sigma$. The regularization term and the gradient can be expressed as
\begin{equation}\label{eq:reg_RED}
\rho_{\bm\sigma, \nu}(\*x) \;\; = \;\; \frac{1}{2\nu}\*x^T \big( \*x - D_{\bm\sigma}(\*x) \big), \;\quad\quad\; \nabla \rho_{\bm\sigma, \nu}(\*x) \;\; = \;\; \frac{1}{\nu} \big( \*x - D_{\bm\sigma}(\*x) \big)
\end{equation}

\noindent if the denoiser $D_{\bm\sigma}(\cdot)$ is locally homogeneous, i.e. $(1 + \epsilon)D_{\bm\sigma}(\*x) = D_{\bm\sigma}\left( (1 + \epsilon)\*x\right), \;\forall\*x, \epsilon>0$ and the its Jacobian is symmetric, $\*J[D_{\bm\sigma}(\*x)] = \*J[D_{\bm\sigma}(\*x)]^T$. When the denoiser does not follows the above conditions, the relationship between $\rho_{\bm\sigma}(\*x)$ and the gradient expression in Eq. (\ref{eq:reg_RED}) does not hold. In the following Section \ref{sec:Denoi_Score}, an interpretation of RED for generic denoisers is summarized which is based on the score-matching by denoising framework \cite{reehorst2019}.

\subsection{Denoising Score Matching}\label{sec:Denoi_Score}

In medical imaging applications, the true prior $p(\*x)$ is not available and therefore it is not possible to compute directly the MMSE estimator $D_{MMSE, \nu}(\cdot)$. Since calculating ${D}_{MMSE,\nu}$ is infeasible in high dimensions, the main idea is to approximate the optimal MMSE estimator by a generic denoiser, $D_{\bm\sigma}(\cdot)$ parametrized over $\bm\sigma$, even for denoisers not locally homogeneous or with non-symmetric Jacobian. 
By utilizing ${p}(\*x,\nu)$ as prior of $\*x$ for the regularization term $\rho(\*x) = -\log{p}(\*x,\nu)$ in (\ref{eq:pr_GLM}), it is necessary to evaluate the expression for $\nabla_{\*x}\rho({\*x}) = \nabla_{\*x}\log{p}(\*x,\nu)$ and $\nabla^2_{\*x}\rho({\*x})$. Instead, using an efficient denoiser on the class $\+D\triangleq\{D_{\bm\sigma}:\bm\sigma\in\bm\Sigma\}$, although an analytic expression for ${p}(\*x,\nu)$ is either not available or complex to handle, from an algorithmic perspective the prior is not needed but only the gradient and its Jacobian. Therefore, we exploit the approach based on approximating the score instead of the prior directly by invoking the connection between denoising and the score, i.e. the gradient of the log-prior, and its Hessian \cite{hyvarinen2005}.  It was shown in \cite{vincent2011} that if we choose a function such that 
\begin{equation}\label{eq:denoising}
\theta_{\bm\sigma}(\*x) \;\; = \;\; \frac{D_{\bm\sigma}(\*x) - \*x}{\nu}, \;\quad\quad\; \frac{\partial\theta_{\bm\sigma}(\*x)}{\partial\*x} \;\; = \;\;  \frac{1}{\nu}\Big(\*J[D_{\bm\sigma}(\*x)] - \*I\Big)
\end{equation} 

\noindent where $\*J[D_{\bm\sigma}(\*x)]$ is the Jacobian matrix of the denoising function of dimension $n\times n$, then for any $D_{\bm\sigma}$ and $\*x$, the score-matching error is connected to the denoiser approximation error as follows
\begin{equation}\label{eq:scorematch}
\Big\| D_{\bm\sigma}(\*x) - {D}_{MMSE,\nu}(\*x) \Big\|^2  = 
\Big\|\theta_{\bm\sigma}(\*x) - \nabla_{\*x}\log{p}(\*x,\nu)\Big\|^2 \;\;
= \;\; 
\left\|\frac{\partial\theta_{\bm\sigma}(\*x)}{\partial\*x} - \nabla^2_{\*x}\log{p}(\*x,\nu)\right\|^2 
\end{equation}

\noindent It is important to stress that the equivalence between Eq. (\ref{eq:denoising}) and (\ref{eq:scorematch}) holds also in the case where the parameter vector of the denoiser $\bm\sigma$ is not optimal. The Jacobian matrix $\*J[D_{\bm\sigma}(\*x)]$ of dimension $n\times n$ is computationally prohibitive to calculate and store at each iteration $t$ of the solver in high dimensions. In Section \ref{sec:dir_der} we propose a way to approximate the Jacobian-vector products required within the proposed iterative sub-sampled Newton solver.

\subsection{Review of the Newton Method}

\noindent In this section we build the analysis of the sub-sampled Newton method by first summarizing the Newton method. Let $g: \R^n \rightarrow \R$ be a closed, convex, and twice-differentiable function that is bounded below. Given the current iterate $\*x^t\in\mathcal{C}$, the Newton method generates a new iterate $\*x^{t+1}$ by performing a constrained minimization of the second-order Taylor expansion 
\begin{equation}\label{eq:New_Taylor}
{\*x}^{t+1} = \arg\min_{\*x\in\mathcal{C}} \left[ \frac{1}{2} \langle\*x - {\*x}^t, \nabla^2_{\*x^t} g({\*x}^t, \*y)(\*x - {\*x}^t) \rangle + \langle \nabla_{\*x^t} g({\*x}^t, \*y), \*x - {\*x}^t \rangle \right] 
\end{equation}

\noindent where for the traditional Newton's method $\nabla^2_{\*x^t} g({\*x}^t, \*y)$ represents the full $n\times n$ Hessian matrix. In the following derivations, we simplify the notation of the gradient and the Hessian with $\nabla g({\*x}^t)$ and $\nabla^2 g({\*x}^t)$ respectively, hiding the dependence over $\*y$. By setting the gradient respect to $\*x^t$ equal to zero, it follows the fixed-point update for the exact full Newton's method
\begin{equation}\label{eq:Newton_updates}
\nabla^2 g\big({\*x}^t\big)\big(\*x - {\*x}^t\big) \;\; = \;\; - \nabla g\big({\*x}^t\big)  \quad
\Rightarrow \quad \*x^{t + 1} 
\;\; = \;\; \*x^t - \alpha^t \left[\nabla^2 g\big({\*x}^t\big)\right]^{-1} \nabla g\big({\*x}^t\big) 
\end{equation}

\noindent The two forms of updates presented in (\ref{eq:Newton_updates}) highlight that 2 different types of solvers can be employed; in the first case, we solve a linear system in order to estimate the vector $\*p^t=\*x - {\*x}^t$ and this is generally obtained by using the CG algorithm. Alternatively, it is possible to directly invert the square $n\times n$ Hessian matrix. The step-size $\alpha^t$ in (\ref{eq:Newton_updates}) can be assumed constant, $\alpha^t = 1, \; \forall t\in\mathbb{N}$.

In the following Section, we describe the proposed Denoising-IHS algorithm which uses a stochastic approximation of the Hessian $\nabla^2 g({\*x}^t)$ in order to compute $\*p^t$ with reduced computational cost, instead of calculating an $n\times n$ Hessian matrix and solving the linear system in (\ref{eq:Newton_updates}) exactly which would be remarkably expensive in terms of time. We assume that the gradient $\nabla g({\*x}^t)$ is computed exactly and the Hessian is approximated by random projection to lower dimension $s<n$ by applying a sketching matrix $\*S^t\in\R^{s\times n}$ such that $\frac{1}{s}\Expect\big[(\*S^t)^T \*S^t\big] = \*I_n$.

\section{Proposed Denoising-based Sub-sampled Newton Method}\label{sec:proposed_alg}

The Denoising-IHS algorithm minimizes the original MAP-driven regularized composite cost function of the form $g(\*x) = f(\*x) + \rho_{\bm\sigma}(\*x)$ in Eq. (\ref{eq:pr_GLM}) and we utilize a partially sketched Newton update which means that a sketch of the Hessian $\nabla^2_{\*x} f(\*x,\*y)$ is performed while retaining the exact form of the Hessian associated with the regularizer, i.e. $\nabla^2\rho_{\bm\sigma}(\*x)$. The reason behind this choice is while it is possible to reduce the computational complexity related to the loss function by using fast non-uniform rows sampling, for generic black-box denoising function the computational cost of applying the sketching to the regularizer is higher than using the original function. The algorithm assumes that the square root of the Hessian matrix of dimension $m\times n$, indicated as $\left[\nabla^2 f({\*x})\right]^{\frac{1}{2}}$, can be computed for the CT model in (\ref{eq:HessianCT}), it follows that 
\label{eq:square_root_H} $\left[\nabla^2 f({\*x}^t)\right]^{\frac{1}{2}} = \bm\Sigma^{-\frac{1}{2}}\*A$

\noindent being $\bm\Sigma$ a diagonal matrix. 
The sketched Hessian of the loss function $f(\*x)$ is constructed as
\begin{equation}\label{eq:G_sketch_H}
\*H^t_f(\*x^t) \;\; = \;\; \left(\*S^t \left[\nabla^2 f({\*x}^t)\right]^{\frac{1}{2}}\right)^T \*S^t \left[\nabla^2 f({\*x}^t)\right]^{\frac{1}{2}} \;\; = \;\; \*G(\*x^t)^T\*G(\*x^t), \quad\quad \*S^t\in\R^{s\times n}
\end{equation}

\noindent with 
$\*G(\*x^t) = \*S^t \left[\nabla^2 f({\*x}^t)\right]^{\frac{1}{2}}$, 
while the full Hessian of the denoising-based regularizer is denoted by $\*H^t_{\rho}(\*x^t) = \nabla^2\rho_{\bm\sigma}({\*x}^t)$ with $\*H^t_{\rho}(\*x^t)\succeq 0$ a positive semi-definite matrix. 
Eq. (\ref{eq:G_sketch_H}) represents an unbiased estimator for the full Hessian of the loss
\begin{equation}
\mathbb{E}_t \left[ \left(\left[\nabla^2 f({\*x}^t)\right]^{\frac{1}{2}}\right)^T (\*S^t)^T \*S^t\left[\nabla^2 f({\*x}^t)\right]^{\frac{1}{2}}\right] = \nabla^2 f({\*x}^t)
\end{equation}

\noindent where the expectation is taken over the iteration index $t$ conditioned over the random sketches selected at previous iterations $[1,\ldots, t-1]$. The splitting leads to the partially sketched update
\begin{eqnarray}\label{eq:NewtonSketch}
{\*x}^{t+1} & = & \arg\min_{\*x\in\mathcal{C}}\left[\frac{1}{2}(\*x - {\*x}^t)^T \*H^t(\*x^t)(\*x - {\*x}^t) + \; \langle\nabla g({\*x}^t), \*x - {\*x}^t \rangle \right] \\
\*H^t(\*x^t) & = & \*H^t_f(\*x^t) + \*H^t_{\rho}(\*x^t) = \big(\*S^t[\nabla^2 f({\*x}^t)]^{\frac{1}{2}}\big)^T\*S^t[\nabla^2 f({\*x}^t)]^{\frac{1}{2}} + \nabla^2\rho_{\bm\sigma}({\*x}^t) \nonumber
\end{eqnarray}

\noindent By substituting into (\ref{eq:NewtonSketch}) the definition for the loss in (\ref{eq:NLPoisson}) and the denoising-based regularizer in (\ref{eq:denoising}), the expressions for the gradient and the Hessian become
\begin{eqnarray}\label{eq:HessianSketch}
\nabla g({\*x}^t) & = & \nabla f({\*x}^t, \*y) + \nabla \rho_{\nu, \bm\sigma}(\*x^t) \;\; = \;\; \*A^T\bm\Sigma^{-1}\left(\*A\*x^t - \*y\right) + \frac{1}{\nu}\big(D_{\bm\sigma}(\*x^t) - \*x^t\big) \nonumber\\
\*G(\*x^t) & = & \*S^t[\nabla^2 f({\*x}^t)]^{\frac{1}{2}} \;\; = \;\; \*S^t\left[\bm\Sigma^{-\frac{1}{2}}\*A\right] \nonumber\\
\*H^t(\*x^t) & = & \*G(\*x^t)^T\*G(\*x^t) + \frac{1}{\nu}\Big(\*J[D_{\bm\sigma}(\*x^t)] - \*I\Big) 
\end{eqnarray}
The update is $\*x^{t+1} = \*x^t + \*p^t$ where $\*p^t$ is obtained as solution of the system of equations
\begin{equation}\label{eq:appr_New}
\*H^t(\*x^t)\*p^t = - \nabla g({\*x}^t)
\end{equation}

\noindent We consider to solve the system (\ref{eq:appr_New}) approximately using the Conjugate gradient method. We focus on the analysis of the accuracy of the solution and the complexity. The inner update is
\begin{equation}
\*x^{t+1} = \*x^t + \*p^t_r 
\end{equation}

\noindent where $\*p^t_r$ is the approximate solution after $r$ inner iterations of the CG solver. We consider the case where the gradient is fully computed which is a realistic scenario in case of separable loss function, and the Hessian is non-uniformly sampled. 
To perform an analysis of the complexity of the Sketch Newton-CG method, we remind that for a symmetric positive definite linear system the convergence of CG is not uniform but depends on the difference between the minimum and maximum eigenvalues. In detail, the convergence rate of the CG method varies at every iteration depending on the spectrum $\{\lambda_1\leq \lambda_2 \leq ,\ldots, \leq \lambda_n\}$ of the positive definite matrix $\*H^t(\*x^t)$. After $r$ steps of the CG method applied to the linear system in (\ref{eq:appr_New}), the iterate $p^t_r$ satisfies
\begin{equation*}\label{eq:CG_rate}
\|\*p^t_r - \*p^t_* \|^2_{\*C} \leq \left( \frac{\lambda_{n-r} - \lambda_1}{\lambda_{n-r} + \lambda_1} \right)^2 \|\*p^t_0 - \*p^t_* \|^2_{\*C}, \quad \*C = \*H^t(\*x^t)
\end{equation*}

\noindent where $\*p^t_*$ indicates the exact solution and $\| \*x \|^2_{\*C} \triangleq \*x^T\*C\*x$. The linear CG algorithm allows to solve the system (\ref{eq:appr_New}) using only matrix/operator-vector products with $\*H^t(\*x^t)$. therefore, we use CG as an Hessian-free inner solver which does not require to store the full $n\times n$ Hessian matrix. It is possible to obtain sketched Hessian-vector products by performing two matrix-vector products with respect to the sketched square-root Hessian. In particular, we compute for each inner iteration $i\in[1,\ldots, r]$ of the CG solver the following variable vector
\begin{equation}
	\*v_i = \*H^t(\*x^t) \*p_i = \left(\*G(\*x^t)^T\*G(\*x^t) + \frac{1}{\nu}\Big(\*J[D_{\bm\sigma}(\*x^t)] - \*I\Big)\right)\*p_i
\end{equation}

\noindent using an auxiliary vector $\*z$ and performing two steps as follows
\begin{equation}
	\*z_i = \*G(\*x^t)\*p_i, \quad\; \*v_i = \*G(\*x^t)^T \*z + \frac{1}{\nu}\*J[D_{\bm\sigma}(\*x^t)]\*p_i - \frac{1}{\nu}\*p_i
\end{equation}

\noindent Furthermore, by using the CG solver the matrix $\*H^t(\*x^t)$ is fixed through the inner iterations (being $t$ the index of the outer iteration as in Algorithm \ref{Algo:SN_PnP} and $i$ the index of the CG solver). In the following, we will not explicitly indicate the index $i$ associated with the vector $\*p$ for notation simplicity. Since $\*H^t(\*x^t)$ depends by the Jacobian of the denoiser, as shown in (\ref{eq:HessianSketch}), it would be unfeasible to use gradient based solvers which requires to calculate the Jacobian at each inner iteration. Instead, with CG only directional derivatives of the denoiser, i.e. $\*J$, are needed to be computed and not the full Jacobian; in Section \ref{sec:dir_der} we present fast methods to estimate accurately the directional derivative of learned CNN-based denoisers.

\begin{algorithm*}[!h]
	\caption{Sub-sampled Newton method with Regularization by Denoising}\label{Algo:SN_PnP}
	\begin{algorithmic}
		\STATE{\textbf{Initialization}: set $t=0$, $\*r^0=\*0$, $\*x^0=\*0$, $\tau_x^0 = 1$}
		\STATE{---- Outer iterations ----}
		\FOR{$t = 1,\ldots, T_{max}$}  
		\STATE{1) Compute the leverage-scores based non-uniform sampling distribution $\{p_i\}_{i=1}^p$} (\ref{th:block_scores})
		\STATE{2) Compute the sketched Hessian of the loss:}
		\STATE{$\quad\; \*G^t(\*x^t) = \*S^t[\nabla^2 f({\*x}^t)]^{\frac{1}{2}} = \*S^t\left[\bm\Sigma^{-\frac{1}{2}}\*A\right]$}
		\STATE{3) Partial sketched Hessian:}
		\STATE{$\quad\; \*H^t(\*x^t) = \left[\*G^t(\*x^t)\right]^T\*G^T(\*x^t) + \frac{1}{\nu}\Big(\*J[D_{\bm\sigma}(\*x^t)] - \*I\Big)$}
		\STATE{---- Inner iterations ----}
		\STATE{4) Solve the liner system using CG:}
		\STATE{$\quad \; \*H^t(\*x^t)\*p^t = - \nabla g({\*x}^t)$}
		\ENDFOR
		\RETURN $\*x^t$
	\end{algorithmic}
\end{algorithm*}

\noindent Although in the worst case, CG will require $n$ iterations to converge (thus requiring the evaluation of $n$ products $\*H^t(\*x^t)\*p^t$), the behavior CG can achieve significant progress in the minimization of the objective function after a reduced number of iterations $r$. The following Theorem outlines that the proposed Denoising-IHS algorithm \ref{Algo:SN_PnP} can achieve quadratic convergence, with constants $C_1, C_2, C_3$ dependent on the minimum and maximum eigenvalues of the full Hessian $\Bar{\*H}(\*x^t))$ and the condition number $\kappa$ of the problem as defined in the proof in Appendix \ref{app:a2}

\begin{theorem}
	Let $\*x^t$ be the iterates generated by the non-uniform sampling Newton-CG algorithm with $r<n$ CG iterations are performed. Then,
	\begin{equation}
	\mathbb{E}_t\| \*x^{t+1} - \*x_* \| \leq  C_1 \| \*x^t - \*x_* \|^2  + \left( {C_2} + C_3  \right) \| \*x^t - \*x_* \|^2
	\end{equation}
\end{theorem}

\section{Denoiser directional derivative estimators}\label{sec:dir_der}

The major computational cost of the inexact Newton solver relies on the calculation at each outer iteration $t$ of the Jacobian matrix, of dimensions $n\times n$, of the denoising function $D_{\bm\sigma}$ in Eq. (\ref{eq:HessianSketch}) whose computation is unfeasible in high dimensions. On a closed look, CG requires instead to calculate explicitly the direction derivative $\*J\left[D_{\bm\sigma}(\*x^t)\right]\*p^t$ that for the $n$-dimensional vector $\*p^t$,can be computed using finite differences at the cost of a single denoising evaluation as in Eq. (\ref{eq:denoising}) through the following expression
\begin{equation}\label{eq:MC_SURE}
\*J\left[D_{\bm\sigma}(\*x^t)\right]\*p^t = \lim_{\epsilon\rightarrow 0}\frac{D_{\bm\sigma}(\*x^t + \epsilon\*p^t) - D_{\bm\sigma}(\*x^t)}{\epsilon} 
\end{equation}

\noindent Although the finite differences approach suffers from  numerical instabilities and also requires the computationally expensive evaluation of non-linear functions, there is an efficient procedure for computing exactly the Jacobian-vector product for neural networks. 

\subsection{Automatic differentiation for CNN denoisers}

The procedure to calculate $\*J[D_{\bm\sigma}(\*x^t)]\*p^t$ is based on deriving the right multiplication of Jacobian matrices. We denote the denoiser $D_{\bm\sigma}(\*x): \R^n \rightarrow \R^m$ and let $\*p\in\R^n$ be the right multiplication vector and $\*u\in\R^n$ be a dummy variable. The left multiplication of the Jacobian matrix, i.e. $\*u^T\*J[D_{\bm\sigma}(\*x^t)]$, can be calculated by back propagation, which is common in reverse mode automatic differentiation packages
\begin{equation}\label{eq:backprop}
\*u^T \frac{\partial D}{\partial \*x^T} = \*u^T \*J\left[D_{\bm\sigma}(\*x^t)\right]
\end{equation}

\noindent We exploit (\ref{eq:backprop}) to compute $\*J[D_{\bm\sigma}(\*x^t)]\*p^t$ by defining $g(\*u) = \*u^T\*J[D_{\bm\sigma}(\*x^t)]$. Since $\*u^T\*J[D_{\bm\sigma}(\*x^t)]$ is a vector in $\R^n$, the mapping of $\*p^t$ can be defined as the function $g(\*u): \R^n \rightarrow \R^n$. We can take derivative of $g(\*u)$ with respect to $\*u$, while providing the left multiplying vector $\*p^t$:
\begin{equation}
\*p^t\frac{\partial g}{\partial \*u} = \*p^t\frac{\partial \left(\*u^T\*J\left[D_{\bm\sigma}(\*x^t)\right]\right)}{\partial \*u} \\ 
= \*p^t \*J\left[D_{\bm\sigma}(\*x^t)\right]^T = \left( \*J\left[D_{\bm\sigma}(\*x^t)\right] \*p^t \right)^T
\end{equation}

\section{Block Leverage Scores Sketching}\label{sec:blk_scores}

We focus on sketches based on non-uniform random block rows sampling which can be easily and efficiently implemented on GPU-based imaging operators such as CT \cite{van2016}. Given a probability distribution
$\{p_i\}_{i=1}^m$, we randomly sample the rows of the matrix $\*A$ a total of $s$ times with replacement from the given probability distribution. The rows of $\*S$ are independent and take the values $\*s^T = \frac{\*e_i}{\sqrt{p_i}}$ with probability $p_i$ for $i=1,\ldots , m$, where $\*e_i\in\R^m$ is the $i$-th canonical basis vector. The discrete probability distribution is chosen according to the ridge leverage scores. 

\begin{definition}{(Ridge Leverage Scores)}
	Given $\*A\in\R^{m\times n}$ and a scalar $\lambda\in\R_+$, then for $i=1,\dots, n$, the $i$-th leverage scores of $A$ is defined as $l^{\lambda}_i(\*A) = \*a_i^T(\*A^T\*A+\lambda\*I)^{\dag}\*a_i$.
\end{definition}

\noindent We consider to estimate the ridge leverage scores of the square root matrix of the Hessian matrix of the loss function, i.e. $\*B = \left[\nabla^2 f({\*x}^t)\right]^{\frac{1}{2}} = \bm\Sigma^{-\frac{1}{2}}\*A$ as defined in (\ref{eq:HessianCT}). 
The idea is to adaptively penalize the leverage scores of the Hessian-related matrix  $\*B$ according to the noise level estimated through the trace of the Jacobian of the denoiser based regularizer $\*H^d_{\rho}$. Defining $\*H_{\rho} = \frac{1}{\nu}\Big(\*J[D_{\bm\sigma}(\*x^t)] - \*I\Big)$, we estimate the trace
\begin{equation}\label{eq:trace_D}
H^d_{\rho} = \mathrm{Trace}[\*H_{\rho}] = \frac{1}{\nu}\Big(\mathrm{Trace}\Big[\*J[D_{\bm\sigma}(\*x^t)]\Big] - 1\Big)
\end{equation}
We exploit the property that $H^d_{\rho} =  \mathrm{Trace}[\*H_{\rho}] = \mathbb{E}[\*n^T \*H_{\rho} \*n]$ if $\*n\sim\mathcal{N}(0, \*I)$ is independent of $\*J$ (where the expectation is taken over $\*n$). Therefore, given a noise realization normally distributed $\*n$, we propose to compute the vector $\*J[D_{\bm\sigma}(\*x^t)]\*n$ instead of the matrix $\*J[D_{\bm\sigma}(\*x^t)]$. The expectation over $\*n$ is calculated using a Monte-Carlo method, i.e. generate $K$ i.i.d. $\mathcal{N}(0, \*I)$ samples vectors, estimate the divergence for each vector and obtain the global divergence by averaging:
\begin{equation}\label{eq:trace}
\mathrm{Trace}\Big[\*J[D_{\bm\sigma}(\*x^t)]\Big] = \frac{1}{m}\sum_{j=1}^K \*n_j^T\left(\frac{D_{\bm\sigma}(\*x^t + \epsilon \*n_j) - D_{\bm\sigma}(\*x^t)}{\epsilon}\right) 
\end{equation}

\noindent Given the vectorized image lying in a high dimensional space, it has been empirically observed \cite{ramani2008} that we can accurately approximate the expected value using only a single random sample, i.e. $K=1$. In all the simulations we have used the Monte Carlo method with $K = 1$.

\begin{definition}{(Block partial leverage scores)}
	Given $\*B\in\R^{N_dN_p\times N_v}$ and the scalar $H_{\rho}^d\in\R_+$ defined in (\ref{eq:trace_D}), let $\{l_i\}_{i=1}^{N_dN_p}$ the leverage scores of the matrix $\*B$. The block partial leverage score for the $i$-th block is defined as 
	\begin{equation}
	l_i^{H^d_{\rho}}(\*B) = \sum_{j = (i-1)N_d + 1}^{iN_d} l_j,\quad i = 1, \ldots, N_p
	\end{equation}
	\noindent and the non-uniform sampling distribution is $\; p_i = \frac{l_i^{H^d_{\rho}}(\*B)}{\sum_{j = 1}^{N_p} l_i^{H^d_{\rho}}(\*B)}, \quad i = 1, \ldots, N_p$
\end{definition}

\begin{theorem}[\cite{xu2016}]\label{th:block_scores}
	Given $\*B\in\R^{N_dN_p\times N_v}$, $H^d_{\rho}\in\R_+$ 
	and $\epsilon\in(0,1)$, let $\{l_i^{H^d_{\rho}}(\*B) \}_{i=1}^{N_p}$ be the block partial leverage scores and $p_i\frac{l_i^{H^d_{\rho}}(\*B)}{\sum_{j = 1}^{N_p} l_i^{H^d_{\rho}}(\*B)}$ 
	Let construct the sketch $\*G=\*S\*B$, as in (\ref{eq:G_sketch_H}), by sampling the $i$-th block of $\*G$ with probability $p_i$ and rescaling it by $\frac{1}{\sqrt{p_i}}$. If
	$s \geq 4 \left(\sum_{i=1}^{N_p} l_i^{H^d_{\rho}}(\*B) \right) \log\frac{4N_v}{\delta\epsilon^2}$
	\noindent with probability at least $1-\delta$, then the following condition holds
	\begin{eqnarray}\label{eq:embedding_sc}
	\left\| \left( \*G^T\*S^T \*S\*G + \*H^d_{\rho} \right) - \left( \*G^T\*G + \*H^d_{\rho} \right) \right\|^2 \leq \epsilon \left\| \*G^T\*G + \*H^d_{\rho} \right\|^2 
	\end{eqnarray}
\end{theorem}

\begin{corollary}
	Let focus on the case of one-row sampling, $d=1$ and $\*H_{\rho}$ a scaled identity matrix, $\*H_{\rho} = m\gamma\*I$. Let define $\*G = \*U\bm\Sigma\*V^T = \sum_{i=1}^m\sigma_i\*u_i\*v_i^T$ as the singular value decomposition. Then, the leverage score is 
	\begin{equation}
	l_i^{\gamma} = \*g_i(\*G^T\*G + m\gamma\*I)^{\dag}\*g_i
	= \sum_{j=1}^n \frac{\sigma_j^2}{\sigma_j^2 + m\gamma} u_{ij}^2, \quad i = 1, \ldots , m 
	\end{equation}
	
	\noindent and therefore, the effective $n$-dimension, $n_{eff}$, is much smaller than $n$
	\begin{equation}
	n_{eff}(\*G) = \sum_{i=1}^m l_i^{\gamma} = \sum_{j=1}^n \frac{\sigma_j^2}{\sigma_j^2 + m\gamma} \ll n
	\end{equation}
	\noindent if $m\gamma$ is larger than the majority of singular values of $\*G^T\*G$.
\end{corollary}

\noindent It is worth noting we have used the notion of SVD to show the property of the dimensionality reduction but using the SVD is not practical for calculating the leverage scores. Instead, we will show how to use fast matrix decomposition through the FFT.

\subsection{Convolutional operators}\label{sec:conv_op}

In principle, computing the leverage scores for a generic matrix $\*A$ has a time complexity $\+O(mn^2)$ equivalent to doing an SVD which is computationally demanding in high dimensions. We focus our attention on convolutional operator which can be used to describe the structure of the physical CT operator. We show that the convolutional structure allows to compute efficiently the leverage scores which determines the non-uniform random sub-sampling of the Hessian.

\begin{definition}{(CT Convolutional operator)}\label{eq:def_Radon}
The continuous 2D X-ray Transform is a convolutional linear operator $\+R: L_2(\mathbb{R}^2)\rightarrow L_2([0, \pi) \times\mathbb{R})$ which computes the line integral of a function in the 2D input space. The Fourier central slice Theorem states that $\+R = \+F^{-1}_{\gamma} \Omega_{\omega^{-1}}\+F$ where $\+F$ is the 2D Fourier transform (FT), $\Omega$ is the coordinate transform operator from Cartesian to polar coordinates, $\omega^{-1} = (\gamma\cos\delta, \gamma\sin\delta)$ and $\+F^{-1}_{\gamma}$ is the inverse 1D FT with respect to $\gamma$. 
\end{definition}

\noindent The output of the linear operator is the sinogram that is a function of $\delta$ and the polar space variable $\rho$. Both $\+R$ and $\+R^T\+R$ are normal-convolutional operators since in the frequency domain
\begin{equation}\label{eq:Gram_continous}
\+R^T\+R \; = \; \+F^H (\Omega_{\omega^{-1}})^T(\+F^{-1}_{\gamma})^H \+F^{-1}_{\gamma} \Omega_{\omega^{-1}}\+F \; \stackrel{(a)}{=} \; \+F^H |\mathrm{det} J_{\omega}|\+F \; \stackrel{(b)}{=} \;  \+F^H\+D\left(\frac{1}{|\rho|}\right)\+F
\end{equation}
\noindent where (a) follows from $\+F_{\gamma}$ being an unitary operator and (b) derives from the back-projection CT filter formulation where $J_{\omega}$ defines the Jacobian of $\omega$ and $\+D\left(\frac{1}{|\rho|}\right)$ is the diagonal polar Fourier space operator. From the continuous to the discrete domain, the Fourier-based Radon transform can be written as \cite{matej2004iterative,o2006fourier}: 
\begin{equation}\label{eq:Radon_op}
\*R = \*F^{-1}_{\gamma} \bm\Omega_{\bm\omega^{-1}}\*F
\end{equation}
where $\*F$ is the 2D unitary discrete Fourier transform operator which takes as input the image $\bm\mu$ of dimensions $I\times I$, with $I = \sqrt{N_v}$. The operator $\bm\Omega_{\bm\omega^{-1}}$ performs a discretized version of the continuous coordinate transform in (\ref{eq:def_Radon}) which outputs a matrix of polar coordinate samples that are equally-spaced along $\rho$ at the discrete locations $\{i\Delta_{\rho}\}$ for $i=-\frac{I}{2},\ldots,\frac{I}{2}-1$. The non-uniform FFT operator $\bm\Omega_{\bm\omega^{-1}}\*F$ takes as input the $I\times I$ input image matrix and output a matrix of dimensions $N_{d}\times N_{p}$; $\*F_{\gamma}$ applies the 1D unitary discrete Fourier transform (DFT) matrix separately to each of the radial lines vectors of dimension $N_{dec}$ and it is defined as the Kronecker product between a 1D DFT matrix $\*F_1$ and the identity matrix, i.e. $\*F_{\gamma} = \*F_1\otimes \*I_{N_{p}}$. Therefore, the final output is a vector of dimensions $J=N_{d}\cdot N_{p}$ where $N_{d}$ is the number of detectors and $N_{p}$ is the number of angles (or number of projections) in agreement with (\ref{eq:linCT}). 

\section{Numerical simulations}\label{sec:num_results}
\begin{figure}[!h]
	\centering
	\includegraphics[width=.45\textwidth]{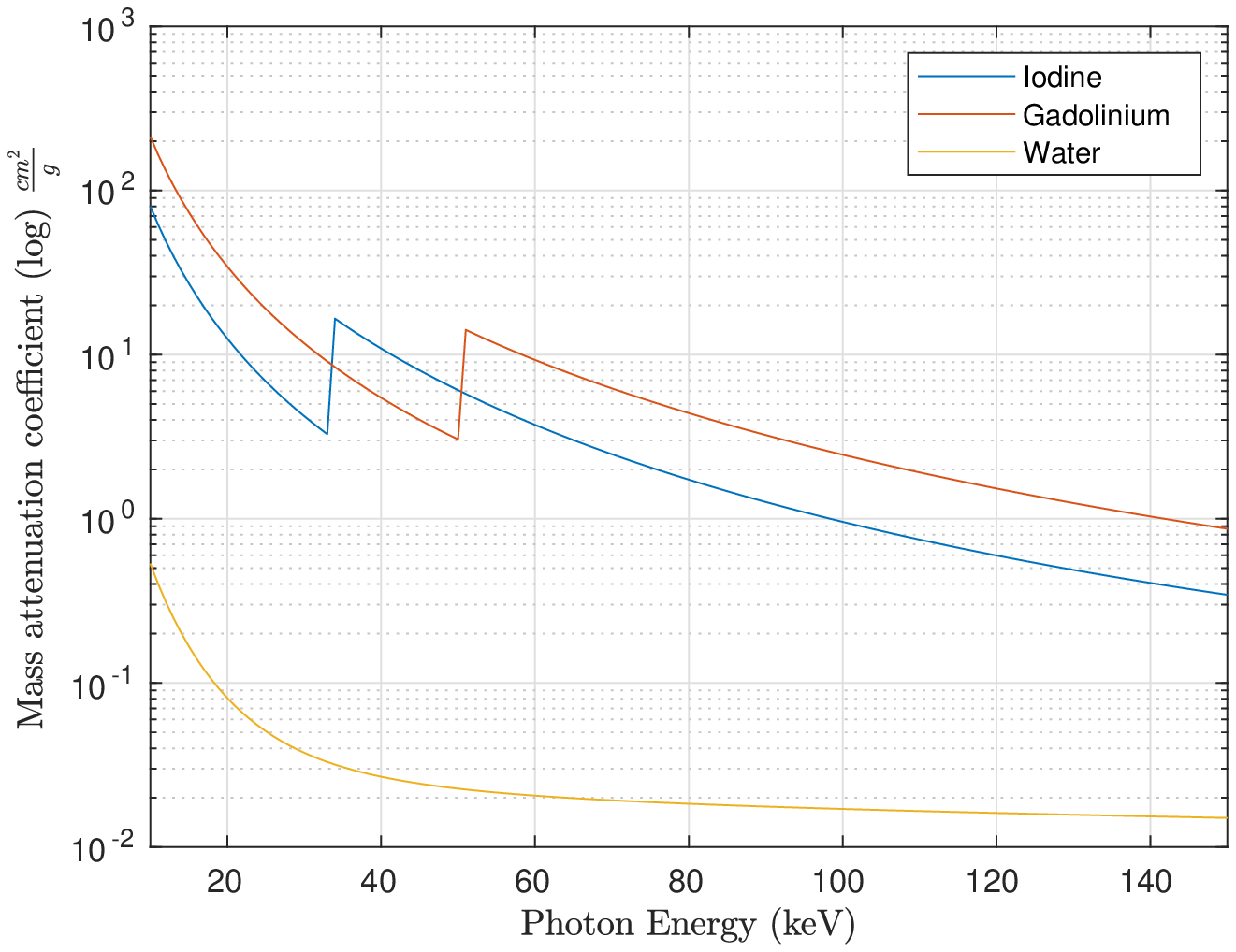}
	\hspace{.0em}
	\includegraphics[width=.45\textwidth]{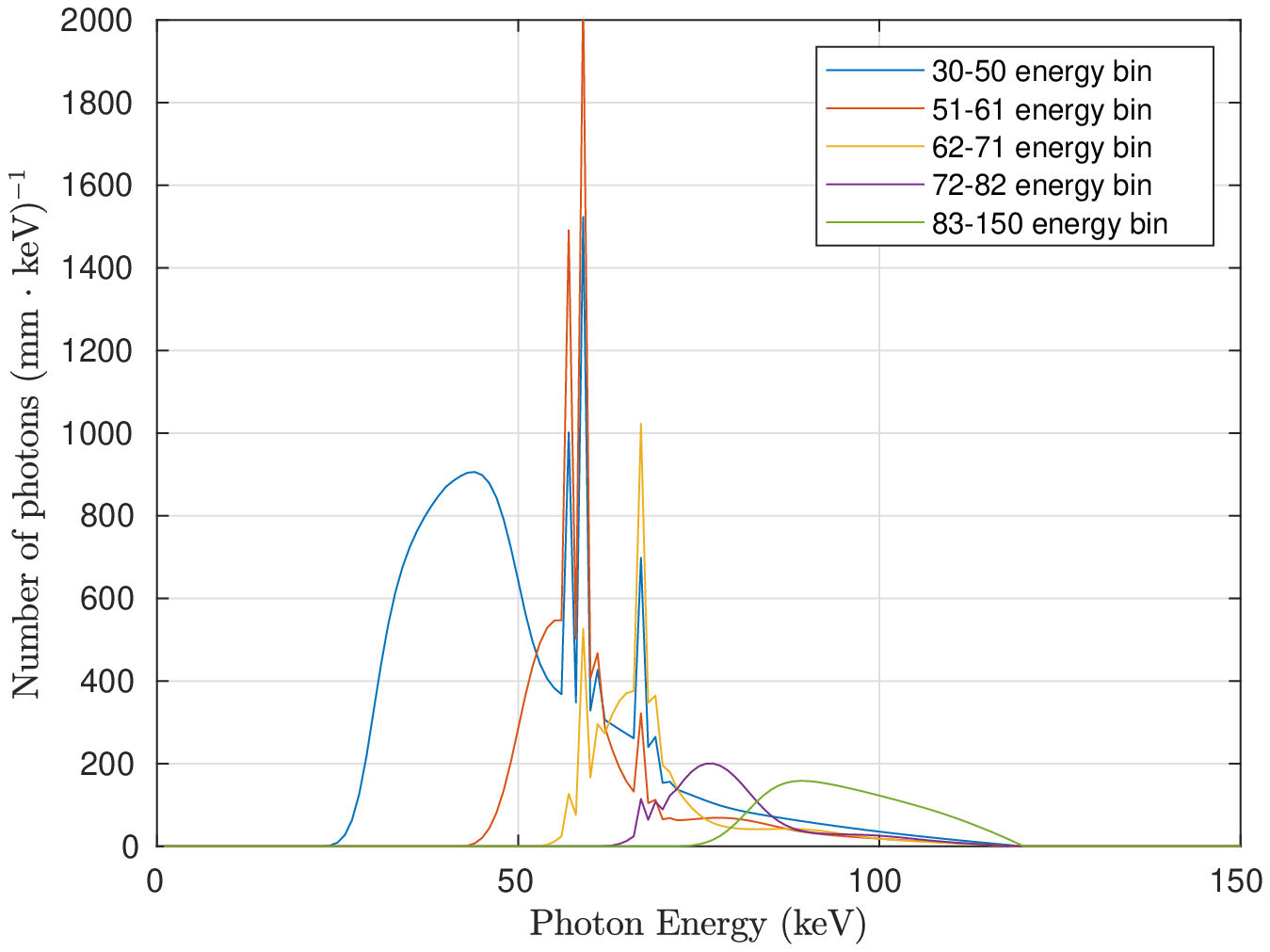}
	\caption{(a) Mass attenuation coefficients of iodine, gadolinium and water as a function of energy and (b) energy spectra of x-ray photons interacting in detector for different energy bins.}\label{fig:energy}
\end{figure}

We have preliminary tested the Denoising-IHS algorithm using the following parameters and geometry; the X-ray spectrum is obtained by using the Spektr \cite{siewerdsen2004} with tube potential $E=80$ keV and initial average photon flux of $I_0=2\cdot 10^3$. Figure \ref{fig:energy} shows the spectrum in the energy interval $[0, 150]$ keV of the five discrete selected energy bins. In Figure \ref{fig:energy}(a), the mass attenuation coefficients of the materials to be estimated are shown; it is interesting to note that the two distinct K-edges associated with the iodine and gadolinium appear at energies corresponding to two distinct energy bins as shown in Figure \ref{fig:energy}(b).  

\begin{figure}[!h]
	\centering
	{\includegraphics[width=.45\textwidth]{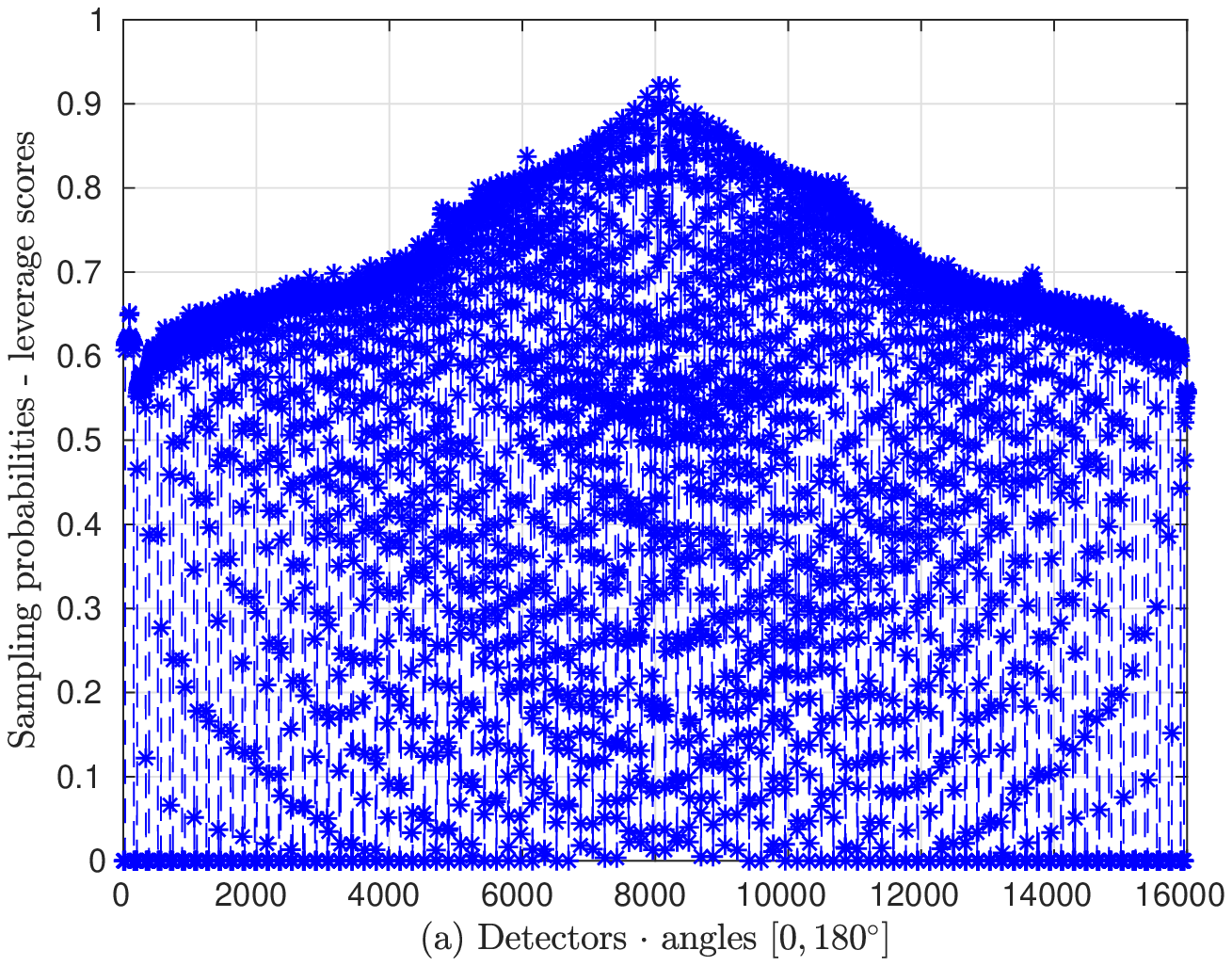}}
	\hspace{.0em}
	{\includegraphics[width=.45\textwidth]{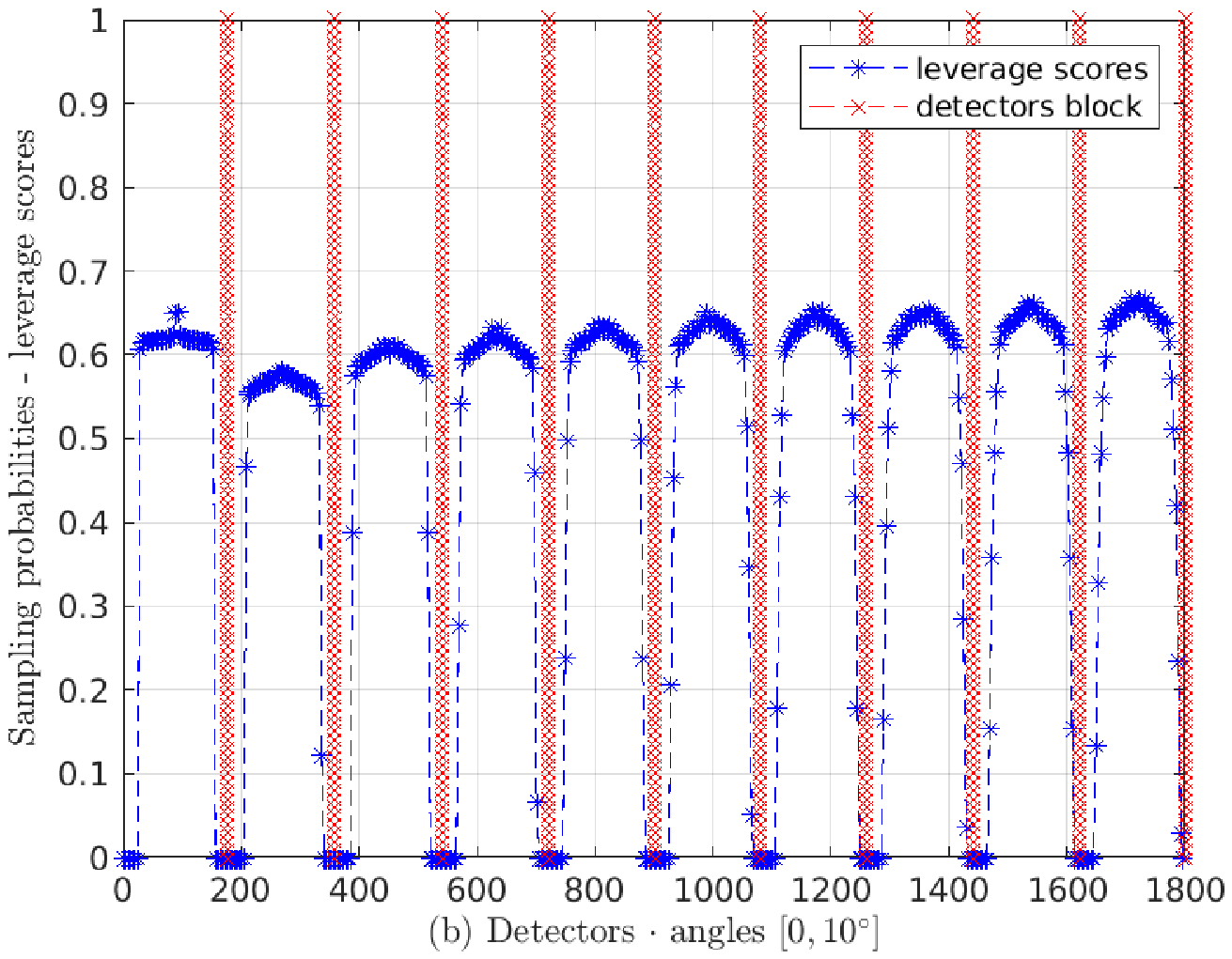}}
	\caption{(a) Leverage scores probabilities on the projection up to $180^{\circ}$ angles and (b) zoom on the structure of the detectors block leverage scores.}\label{fig:lev_scores}
\end{figure}

To generate the multichannel projections, we use the ASTRA-toolbox \cite{van2016} to implement the forward operator associated with a fan-beam CT geometry; the distance from the source to the detector is 45 cm, the source to object distance is 31 cm and the detector cell size is 50 $\mu m$. The pixel resolution of the reconstructed image is $256\times 256$ and the data is simulated on a double grid of $512\times 512$ to avoid the inverse crime, the number of projection angles is set to 360 and the Poisson distributed noise has been added to the projection data. 
Figure \ref{fig:lev_scores}(a) shows the plot of the probabilities of sampling each projection according to the leverage score metric. This plot is obtained by calculating exactly the leverage score of the ASTRA operator offline using the power method. Many projections have a very low score and this shows that a random one-by-one sampling would lose an important quantity of information or it would require to run the algorithm on a number of samples of the order of the projection's dimension.

\begin{figure}[!h]
	\centering
	{\includegraphics[width=.25\textwidth]{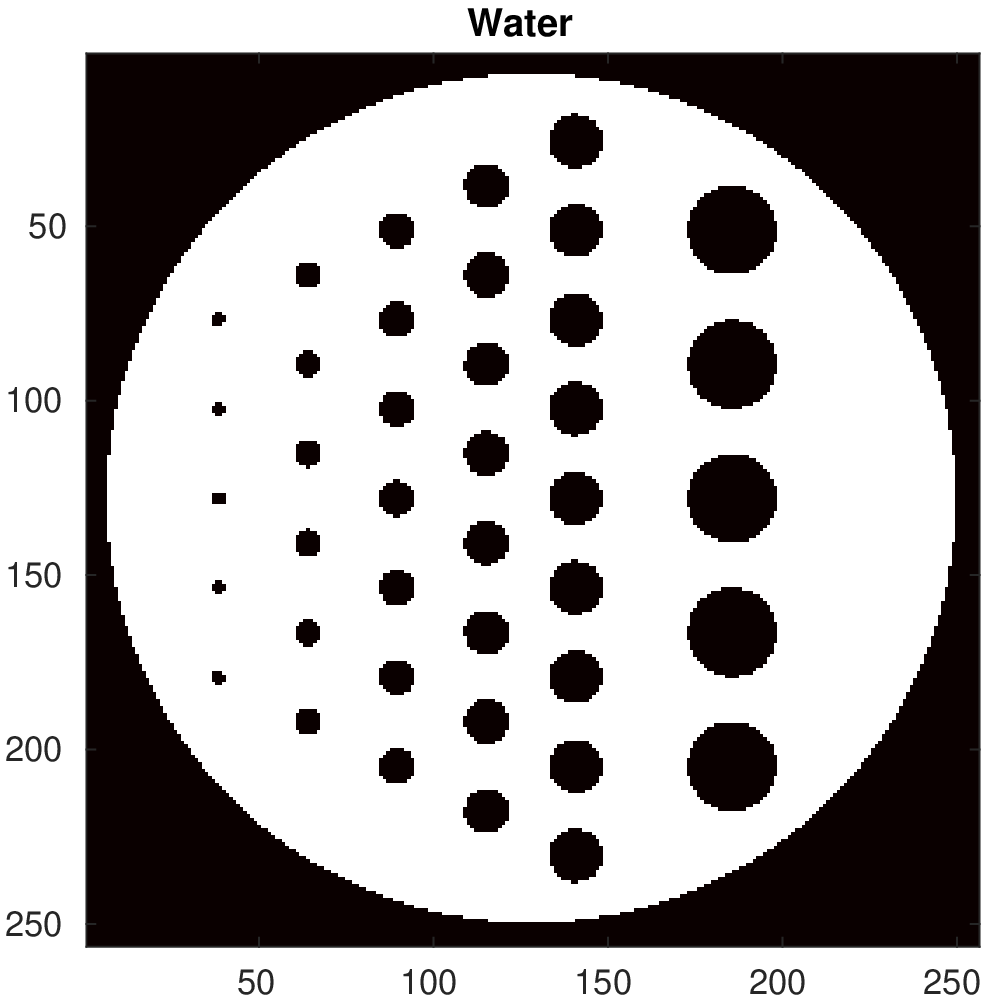}}
	\hspace{.0em}
	{\includegraphics[width=.25\textwidth]{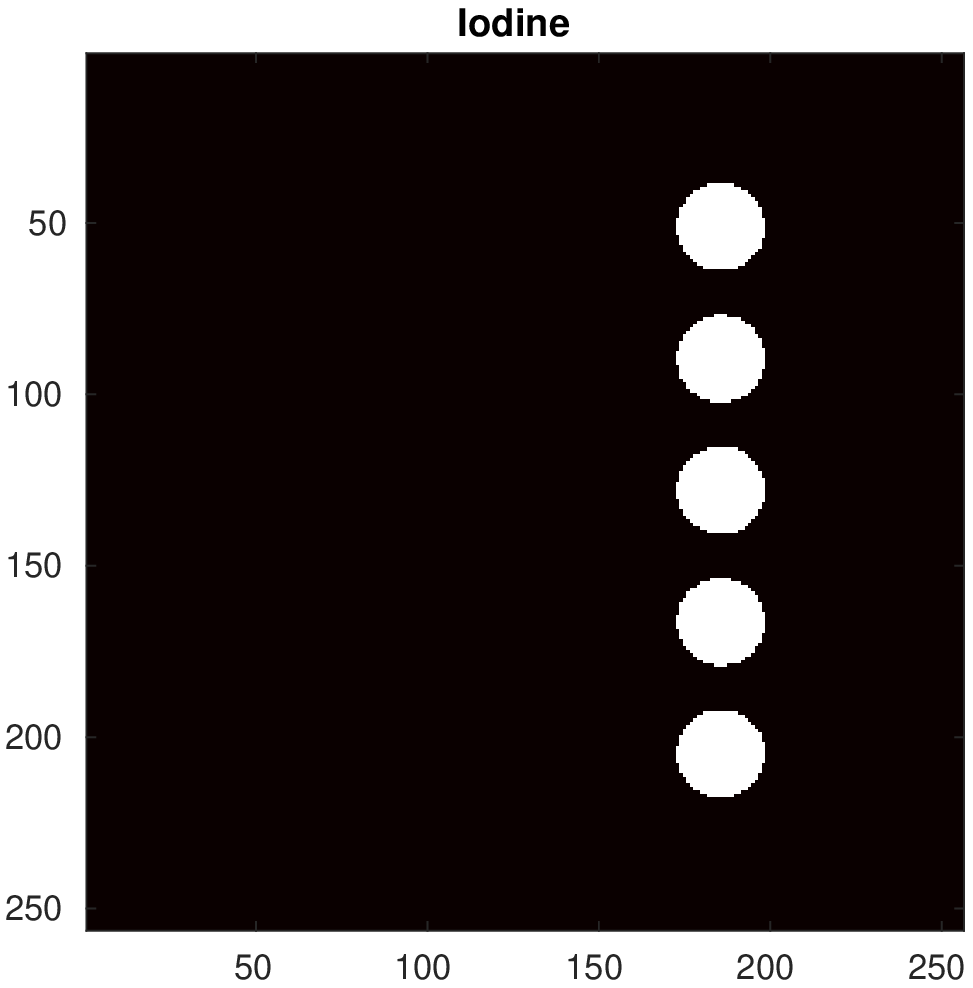}}
	\hspace{.0em}
	{\includegraphics[width=.25\textwidth]{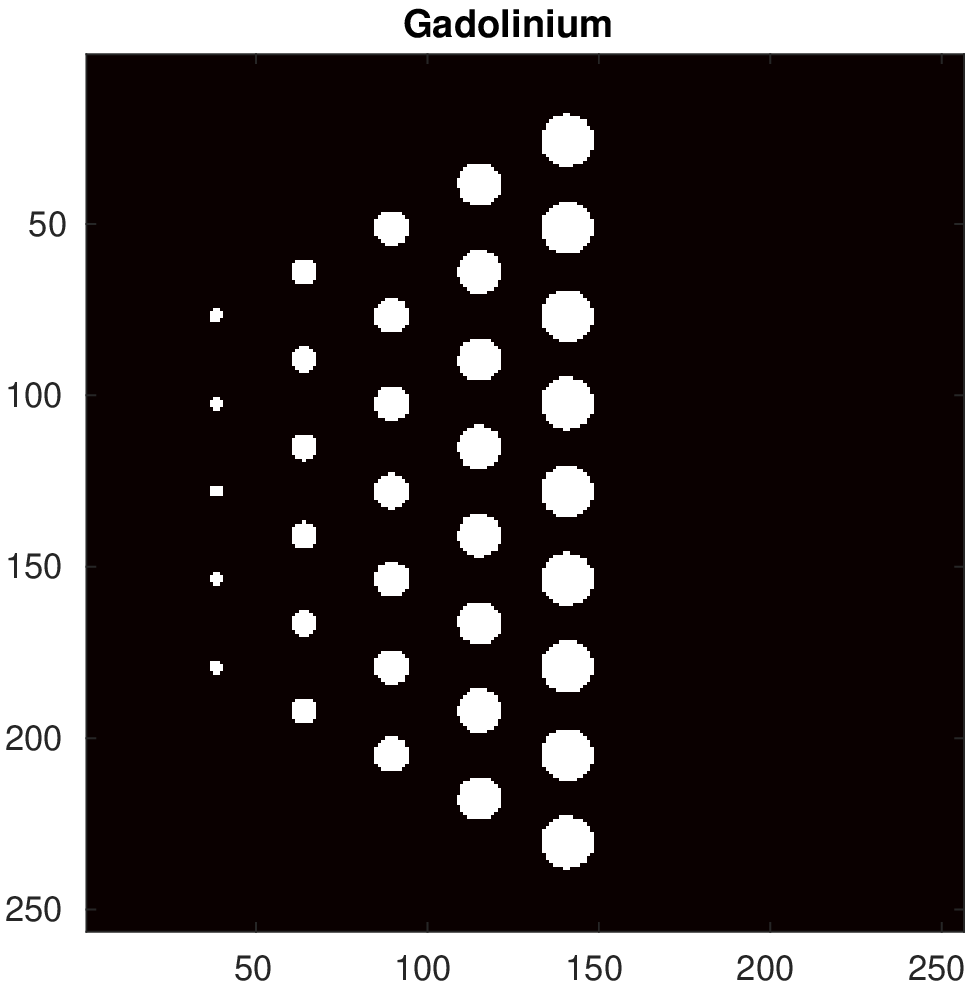}}
	\hspace{.0em}
	{\includegraphics[width=.25\textwidth]{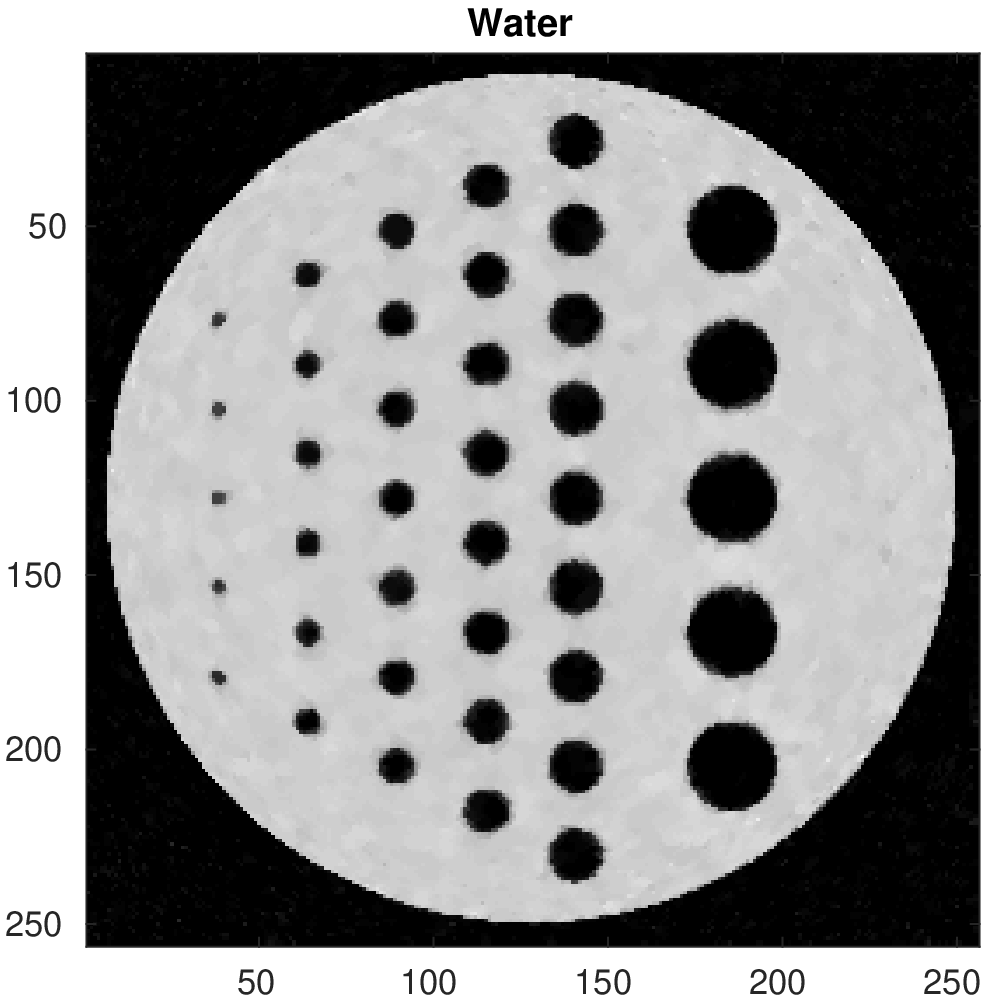}}
	\hspace{.0em}
	{\includegraphics[width=.25\textwidth]{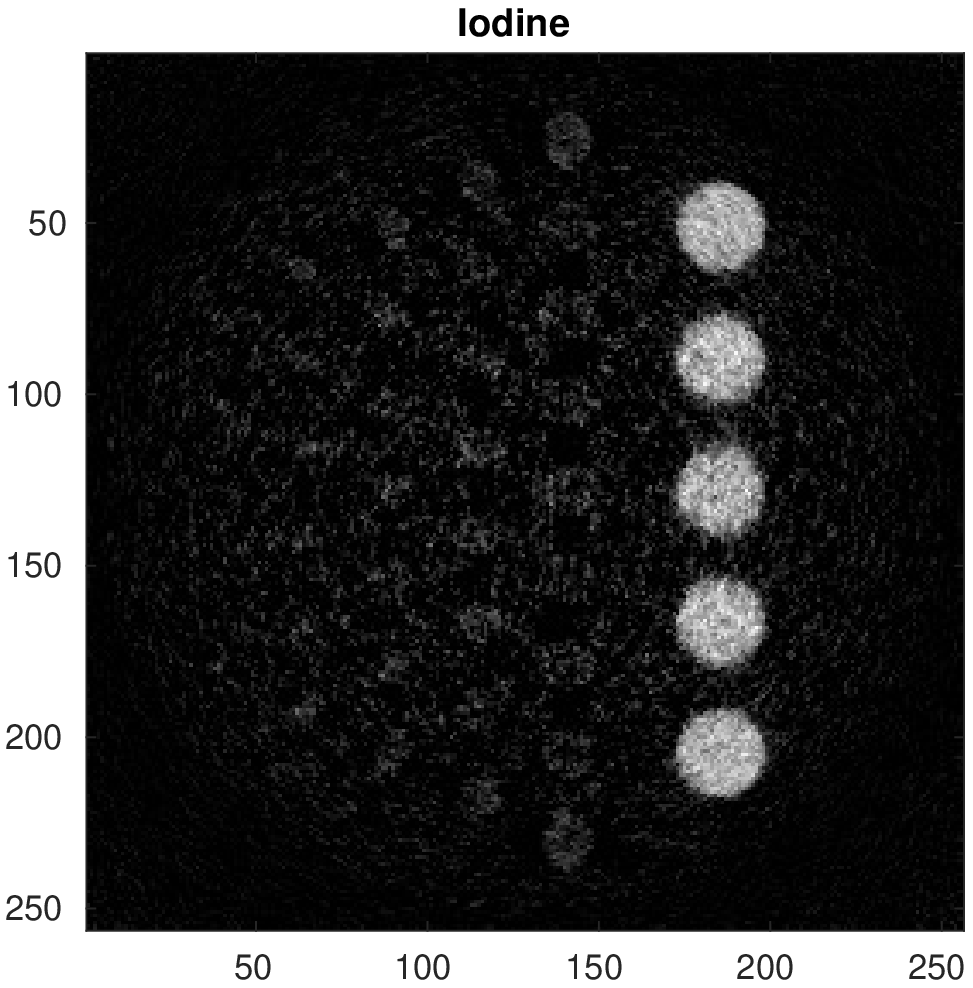}}
	\hspace{.0em}
	{\includegraphics[width=.25\textwidth]{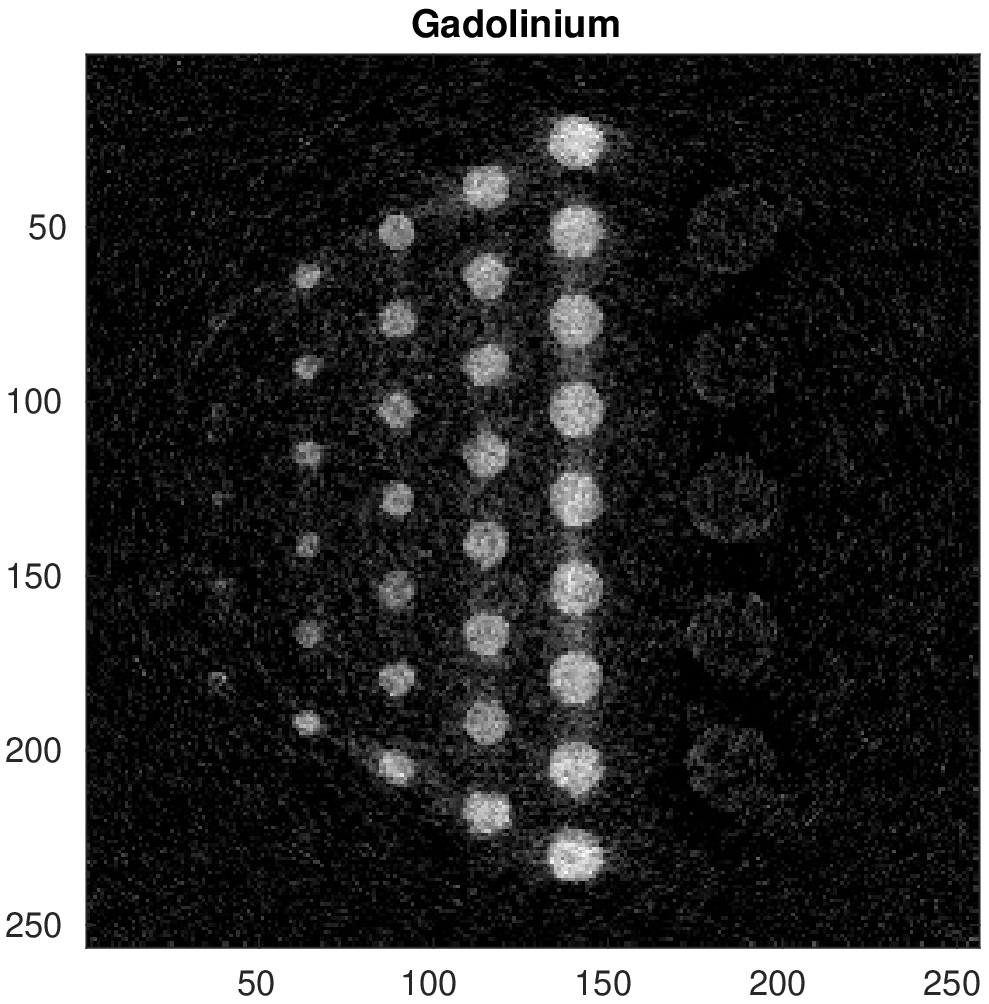}}
	\hspace{.0em}
	{\includegraphics[width=.25\textwidth]{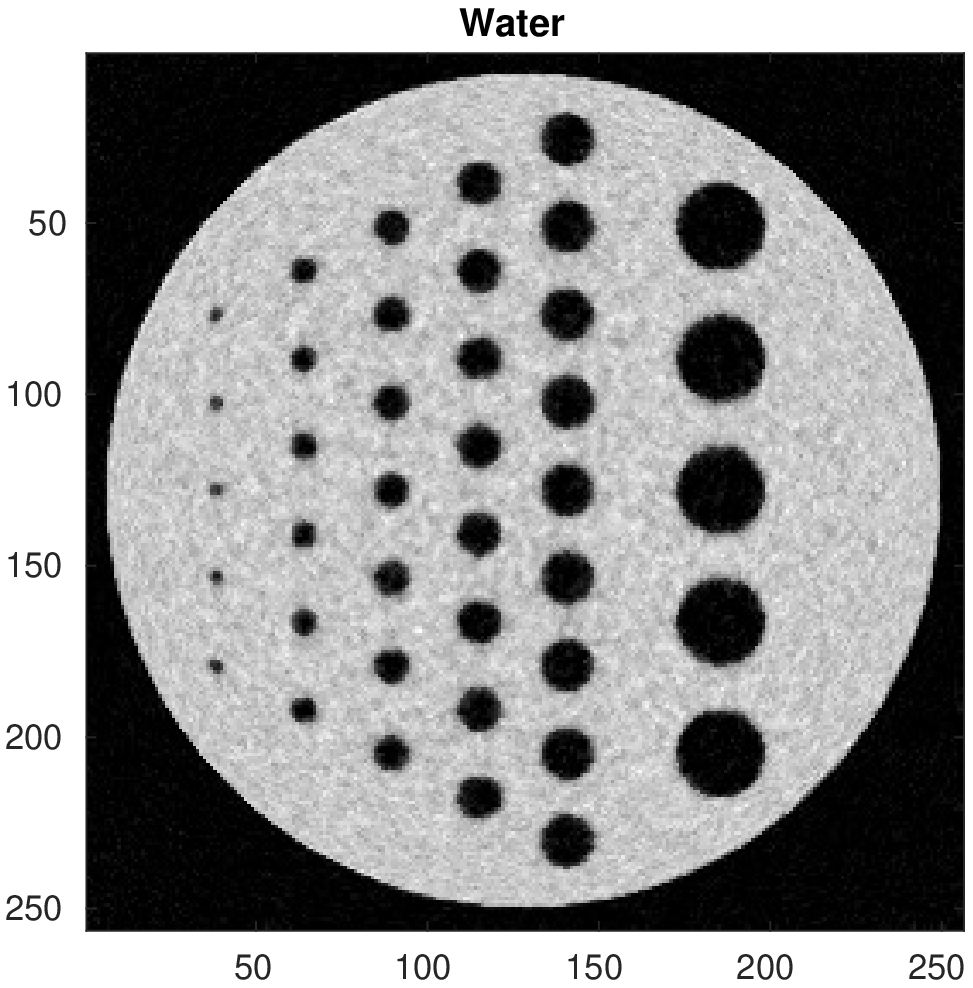}}
	\hspace{.0em}
	{\includegraphics[width=.25\textwidth]{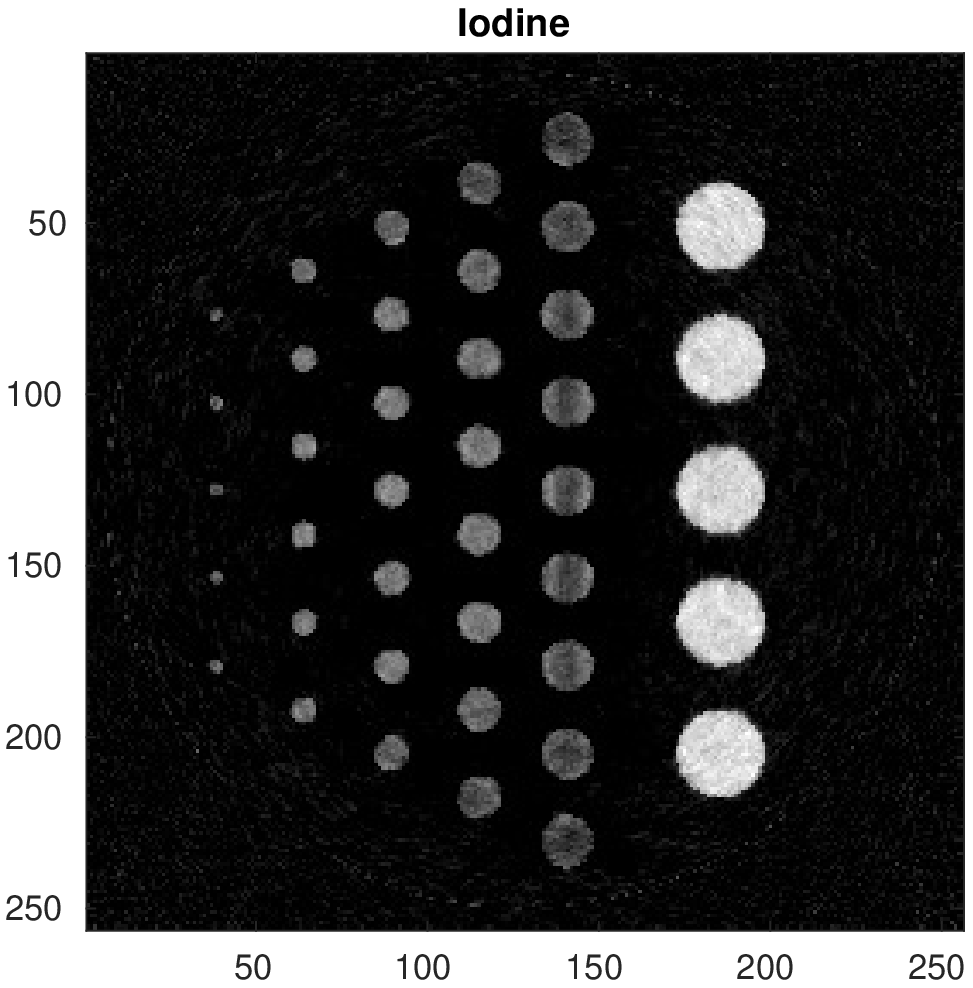}}
	\hspace{.0em}
	{\includegraphics[width=.25\textwidth]{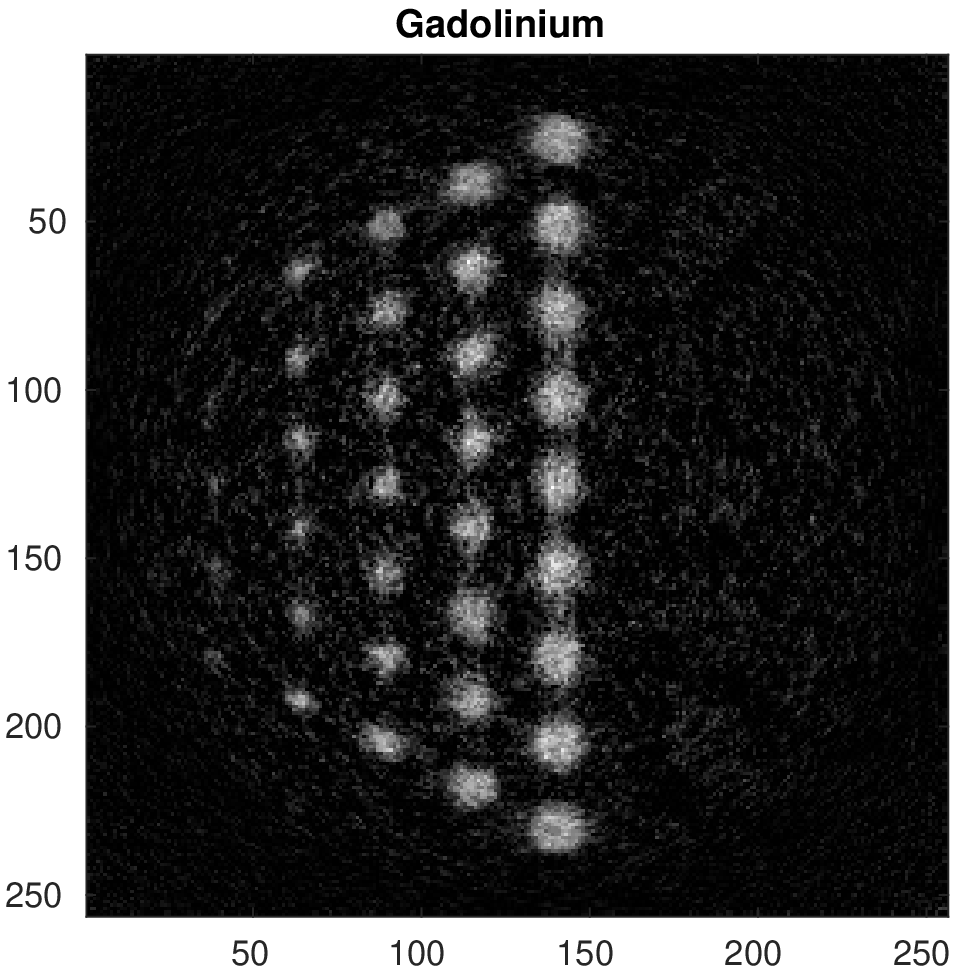}}
	\caption{(a) Ground truth phantom with 3 materials: water, iodine, gadolinium, (b) qualitative materials estimation with the proposed Denoising-IHS algorithm, (c) qualitative materials estimation with algorithm \cite{Long2014} with non-adaptive edge-preserving hyperbola regularizer.}\label{fig:rec_Long}
\end{figure}

\noindent Instead, the block sampling procedure takes advantage, by averaging, of this block structure of the leverage scores induced by the CT operator as depicted in Figure \ref{fig:lev_scores}(b) where each vertical red line corresponds to an adjacent angle and the inner point coincides with the detector array elements. Although this plot represents the exact sampling probabilities, for computational constraint the algorithm estimates these probabilities by exploiting the spectral decomposition of the ASTRA operator through (\ref{eq:Radon_op}). 

We have tested the proposed algorithm on a synthetic phantom composed of three basis materials, water, iodine and gadolinium. The geometric shape of the 2D phantom is build used the model number 4 of the TomoPhantom Matlab software \cite{kazantsev2018Tomo} containing three materials, water, iodine and gadolinium. The phantom is constituted of different circles with increasing diameter as shown in Figure \ref{fig:rec_Long}(a). As denoiser-based regularizer $D_{\bm\sigma}$, we have used the implementation of the U-Net \cite{ronneberger2015} trained on a large set of 1000 images obtained by rotating and placing in different position the circles that are included in the test image. Furthermore, we have used the Autograd software in pytorch to compute the Jacobian-vector product of the U-Net denoiser \cite{paszke2017}. The Hessian block sub-sampling for the Denoising-IHS was set to $s=\frac{m}{3}$.
Figure \ref{fig:rec_Long}(b) shows the quantitative reconstruction of the three materials while Figure \ref{fig:rec_Long}(c) reports the results obtained using the multi-material framework (OS-PWSQS) proposed in \cite{Long2014} which considers a model-based optimization with an edge-preserving hyperbola regularizer. By comparing the images for each materials, it is possible to clearly note how the the Denoising-IHS approach is able to separate the each material and to reduces the noise in the reconstruction, as can be seen in particular for the water and iodine. The iodine image in figure \ref{fig:rec_Long}(c) indicates that algorithm \cite{Long2014} does not manage to fully separate the iodine and gadolinium materials.

\begin{figure}[!h]
	\centering
	{\includegraphics[width=.49\textwidth]{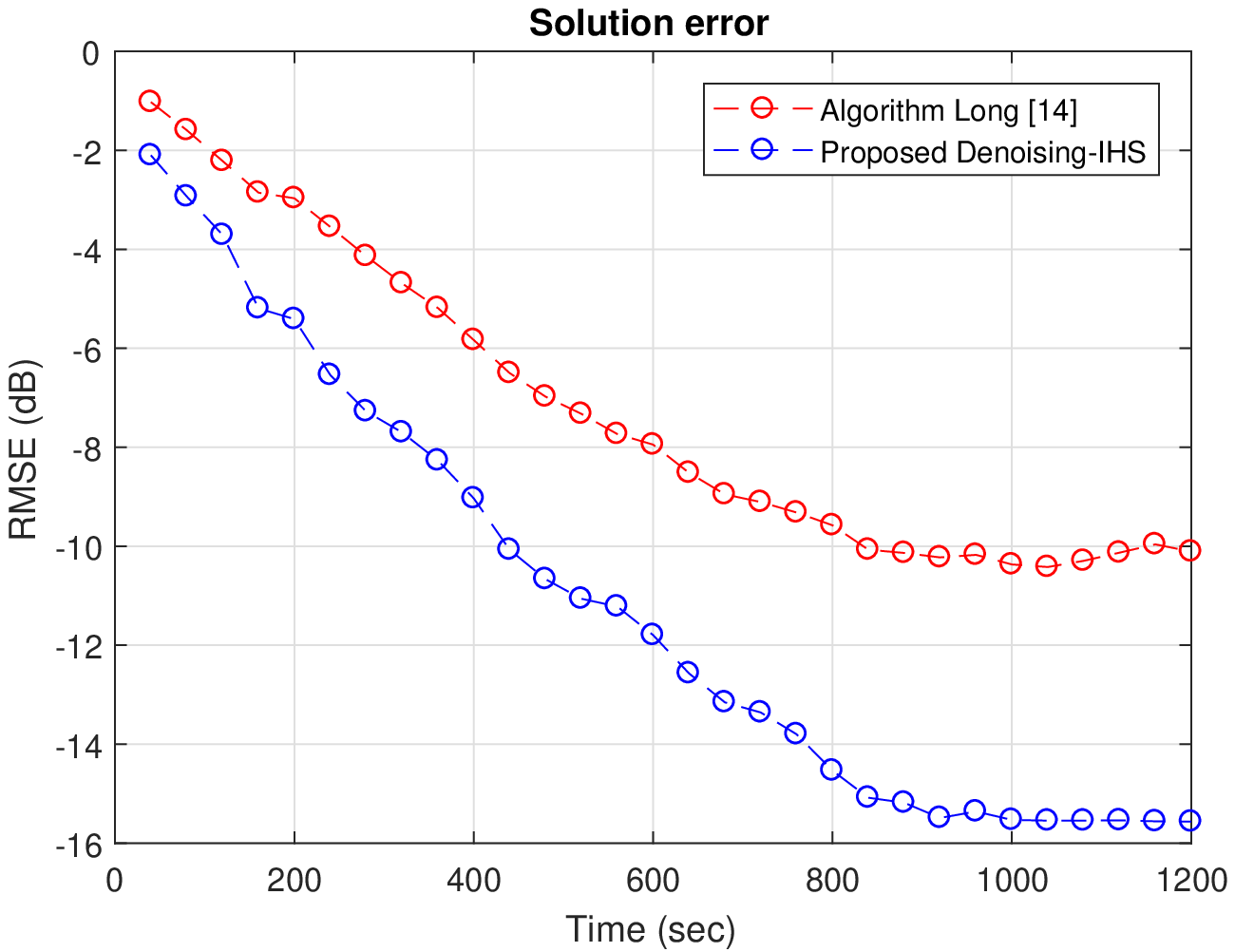}}
	\hspace{.0em}
	{\includegraphics[width=.49\textwidth]{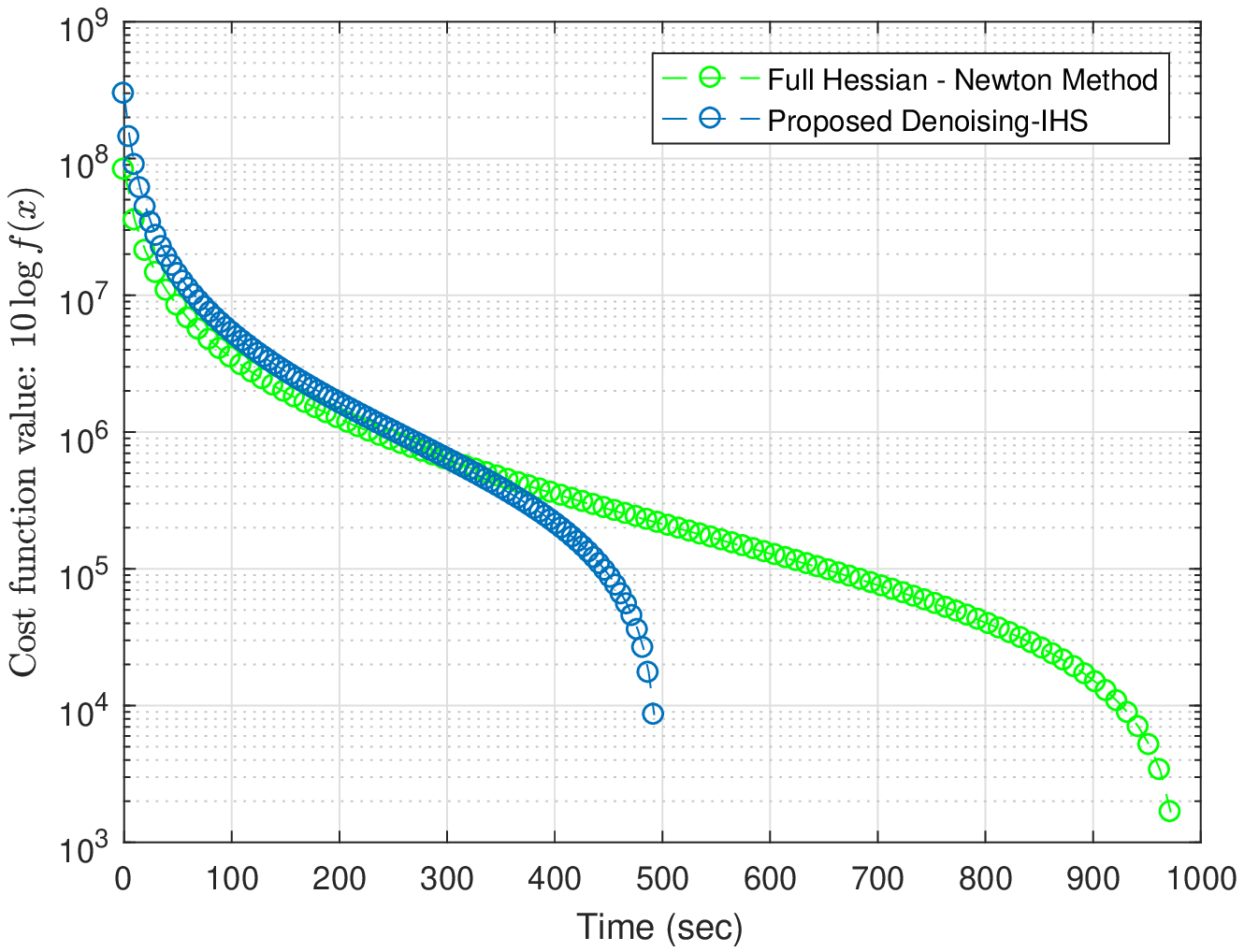}}
	\caption{(a) Relative solution error comparison between the algorithm \cite{Long2014} and the proposed Denoising-IHS algorithm and (b) relative solution error comparison between full Hessian and the proposed Denoising-IHS algorithm.}\label{fig:error_x}
\end{figure}

Figure \ref{fig:error_x}(a) shows the comparison between the root-mean-square error (RMSE), i.e. $\sqrt{\frac{1}{N}(\*x - \*x^*)}$ of the vector $\*x$ containing the concatenation of the estimated materials respect to the ground truth $\*x^*$ between the OS-PWSQS algorithm and the proposed algorithm. The solution achieved by the two algorithms is different because the cost function, both data fidelity term and regularizer, that OS-PWSQS \cite{Long2014} aims at minimizing is different from (\ref{eq:pr_GLM}). By analyzing Figure \ref{fig:error_x}(a), it is possible to highlight that the Denoising-IHS algorithm achieves a more accurate solution. Looking in detail at the blue curve, it is possible to note some local plateaus with sudden reduction of the error and this is in agreement with the fact that in the outer loop the divergence of the denoiser is calculated to update the ridge leverage scores (\ref{eq:trace_D}) which has a cost equivalent of applying the denoiser as stated in (\ref{eq:trace}). 
Furthermore, Figure \ref{fig:error_x}(b) shows the comparison in objective function between the convergence plot (semi-log) in dB of the Denoising-IHS algorithm and the solver which consider the full Hessian of the data loss at each iteration. This plot highlights how the strategy for block sub-sampling of the Hessian leads to a noticeable reduction in computation (almost 2.2 times) while both reach a similar minimum value in the cost function. We can analyze that at early iterations when the estimate is inaccurate both algorithms have similar accuracy/computation trade-off, while going further the Newton solver with full Hessian requires higher computation per iteration and overall it converges close to the Denoising-IHS solution but with increased computational time. 
From Figure \ref{fig:error_x}(b) the time reduction is slightly less compared to the measurement sub-sampling rate and this could be anticipated since the ASTRA operator is not exactly linear in the number of measurements and the Denoising-IHS requires an additional denoising step at each outer iteration for the sampling probabilities update.
  
\begin{figure}[!h]
	\centering
	{\includegraphics[width=.45\textwidth]{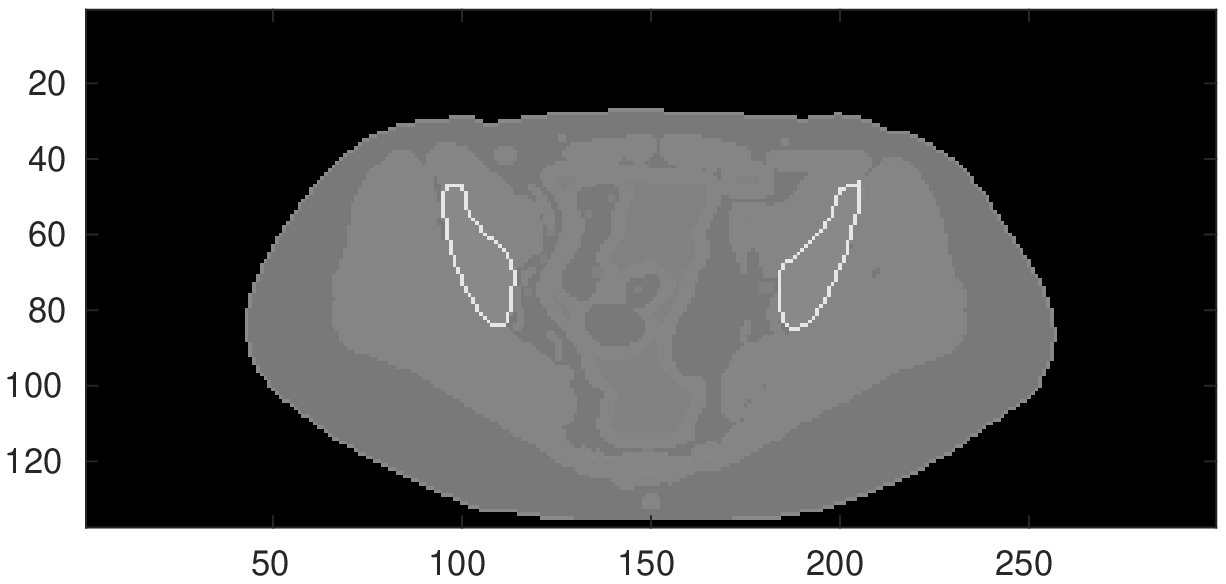}}
	\hspace{.0em}
	{\includegraphics[width=.45\textwidth]{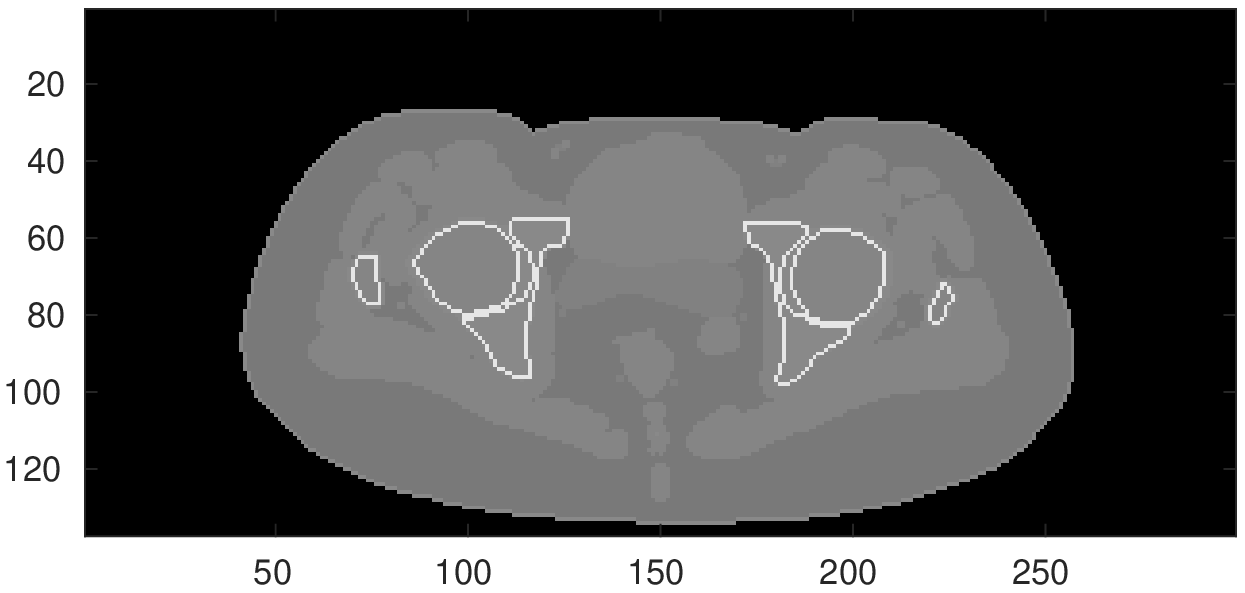}}
	\hspace{.0em}
	{\includegraphics[width=.45\textwidth]{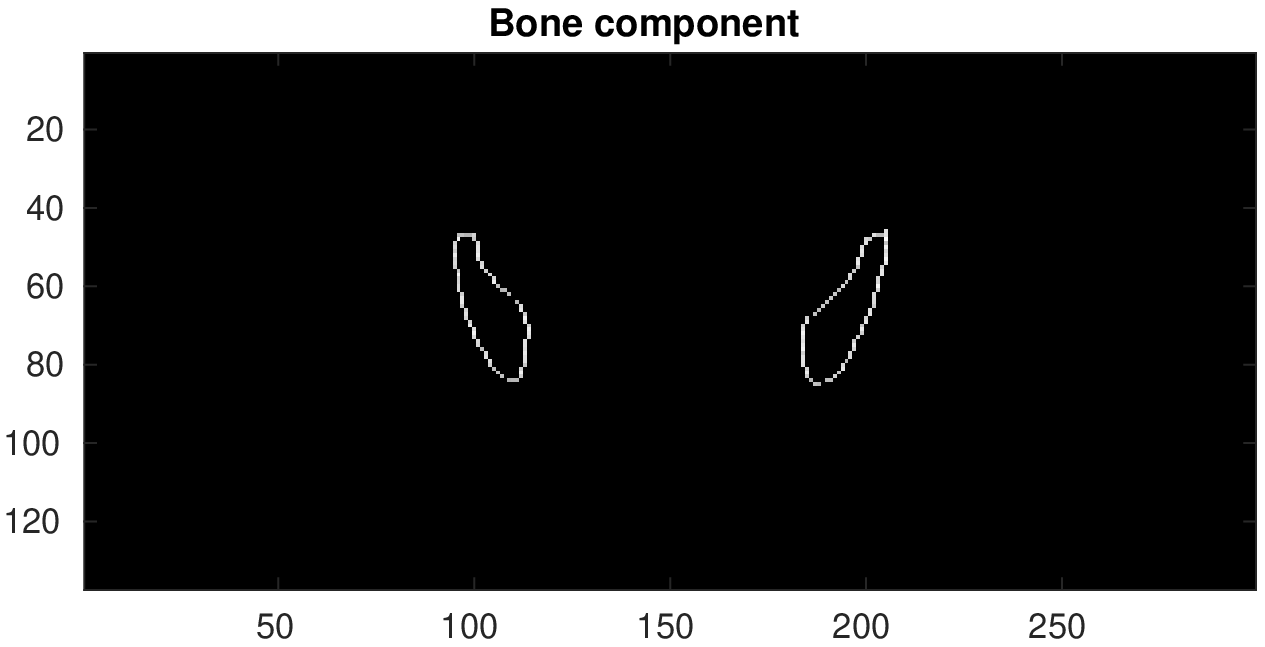}}
	\hspace{.0em}
	{\includegraphics[width=.45\textwidth]{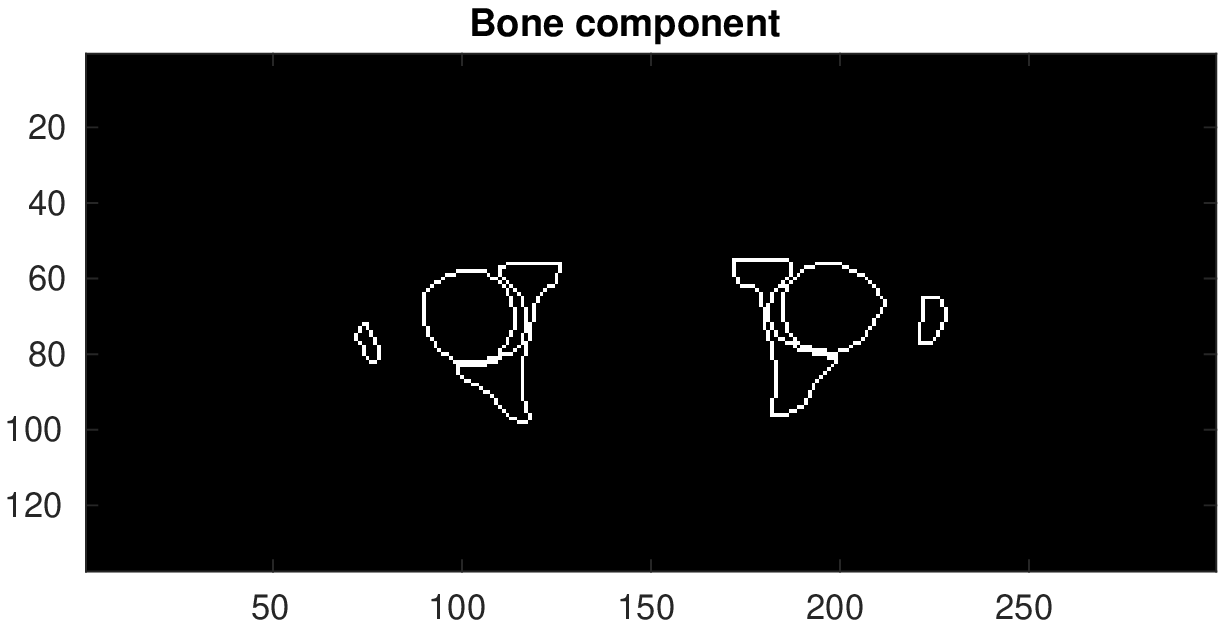}}
	\hspace{.0em}
	{\includegraphics[width=.45\textwidth]{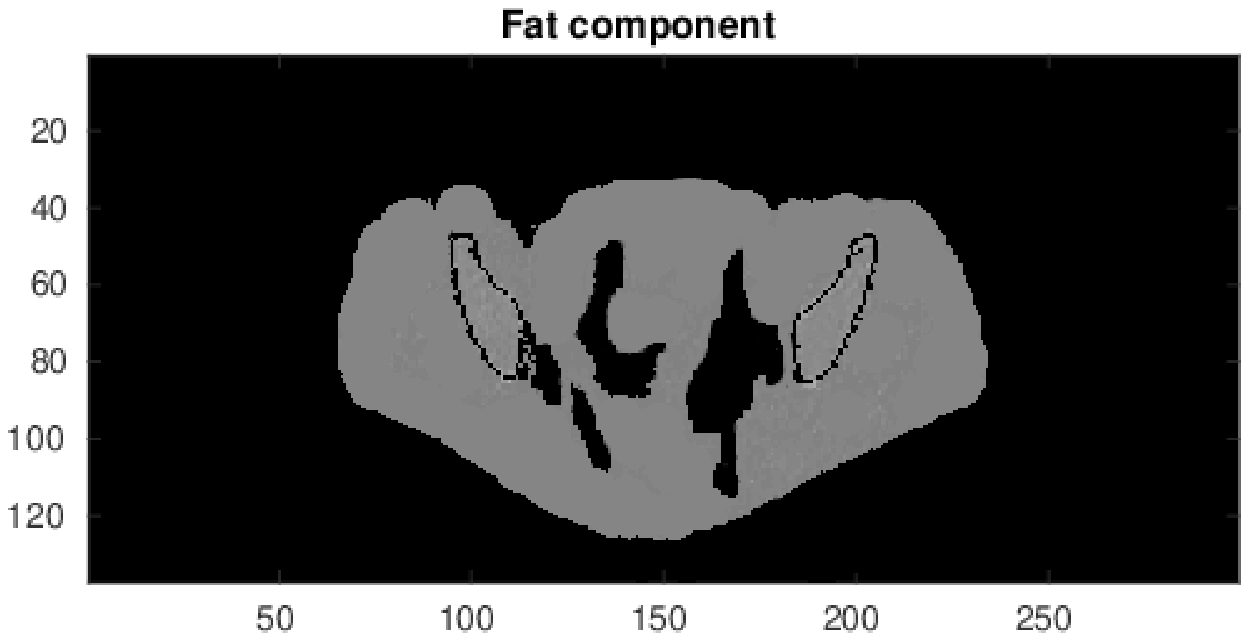}}
	\hspace{.0em}
	{\includegraphics[width=.45\textwidth]{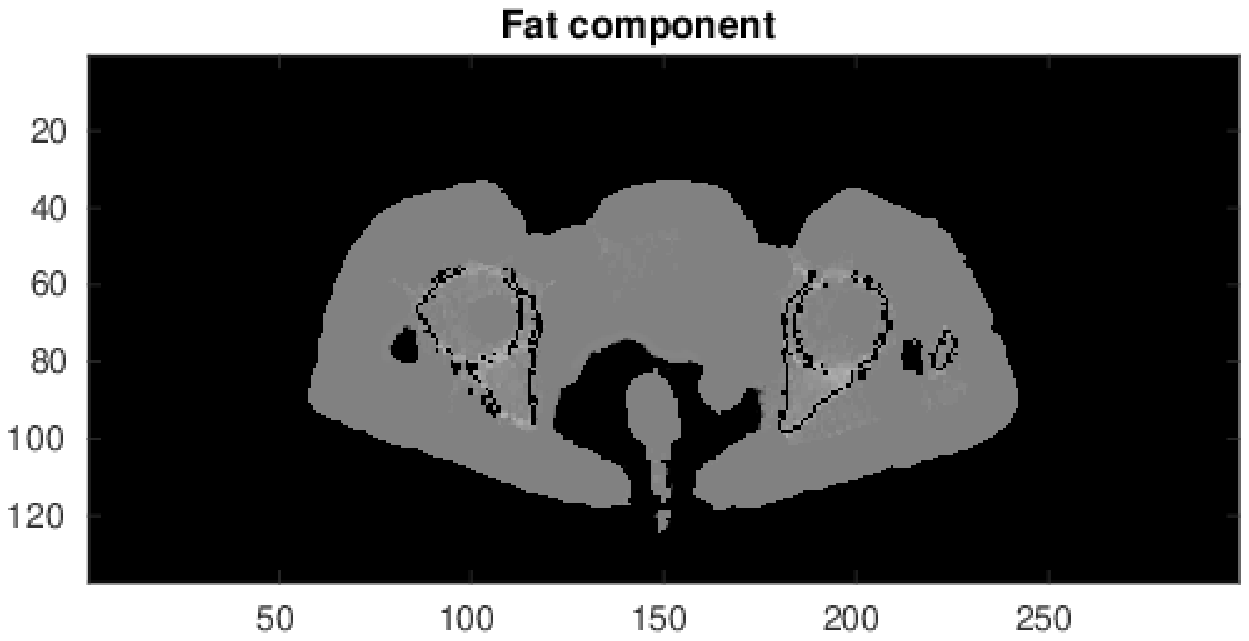}}
	\caption{(a) Testing slices of the pelvis both with a display window $[0, 2]$, (b) reconstruction of the bone material of the two testing images with Denoising-IHS, (c) estimation of the fat material for the two testing images with Denoising-IHS with display window $[0, 0.15]$.}\label{fig:pelvis_fat}
\end{figure}

The following numerical test considers a more structured dataset, derived from the Adult Reference Computational Phantom (ICRP Publication 110, 2009), which is a segmented image of defined chemical composition to represent real tissues \cite{mason2017}. Figure \ref{fig:pelvis_fat}(a) shows two slices selected for testing through the pelvis with image resolution $137\times 299$ with a display window of $[0, 2]$.

\noindent Figure \ref{fig:pelvis_fat}(b) reports the qualitative estimation of the bone component in the testing images; it is possible to note how the Denoising-IHS algorithm is able to tackle the contours of the bone precisely in both images. Figure \ref{fig:pelvis_fat}(c) shows the qualitative estimation of the fat component in the testing images; in this case the Denoising-IHS algorithm manages to estimate correctly most of the complicated structure in the pelvis with some spurious pixels missing.

\begin{figure}[!h]
	\centering
	\includegraphics[width=.6\textwidth]{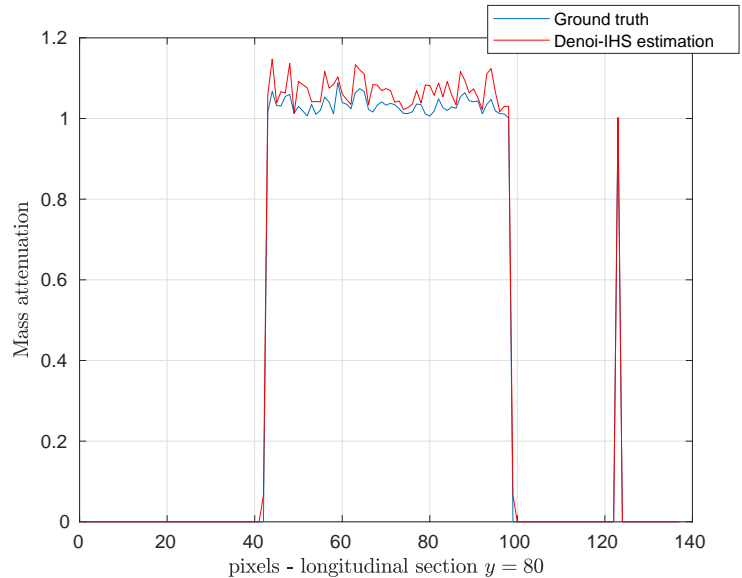}
	\captionof{figure}{Ground truth values of the fat component for the longitudinal section of pixel 80 and the Denoising-IHS estimate.}\label{fig:error_q}
\end{figure}

\noindent To evaluate the quantitative estimation, Figure \ref{fig:error_q} reports the normalized oracle mass attenuation of the fat and the estimated one with Denoising-IHS related to longitudinal section at pixel $y=80$ of the first test image. 

\noindent This plot confirms that the proposed algorithm is able to accurately estimate the mass attenuation with and error which is less than 3\%. In particular, Table \ref{tab:quant} reports the root mean square error (RMSE) for different materials compared to the oracle ground truth. 

\begin{table}[!h]
	\centering
	\begin{tabular}{|l||c|c|c|}
		\hline
		\textbf{Algorithms}				&	\textbf{Bone}	 & \textbf{Fat} & \textbf{Air}\\
		\hline
		\hline
		Denoising-IHS			&	0.031\%	&  0.053\% & 0.025\%		\\
		\hline
		OS-PWSQS			&	0.041\%	 & 0.082\%   & 0.035\%\\
		\hline
	\end{tabular}
	\captionof{table}{RMSE comparison of the material images using Denoising-IHS and OS-PWSQS.}\label{tab:quant}
\end{table}

\section{Preliminary Experimental Results}\label{sec:exp_res}

\begin{figure}[!h]
	\centering
	{\includegraphics[width=.49\textwidth]{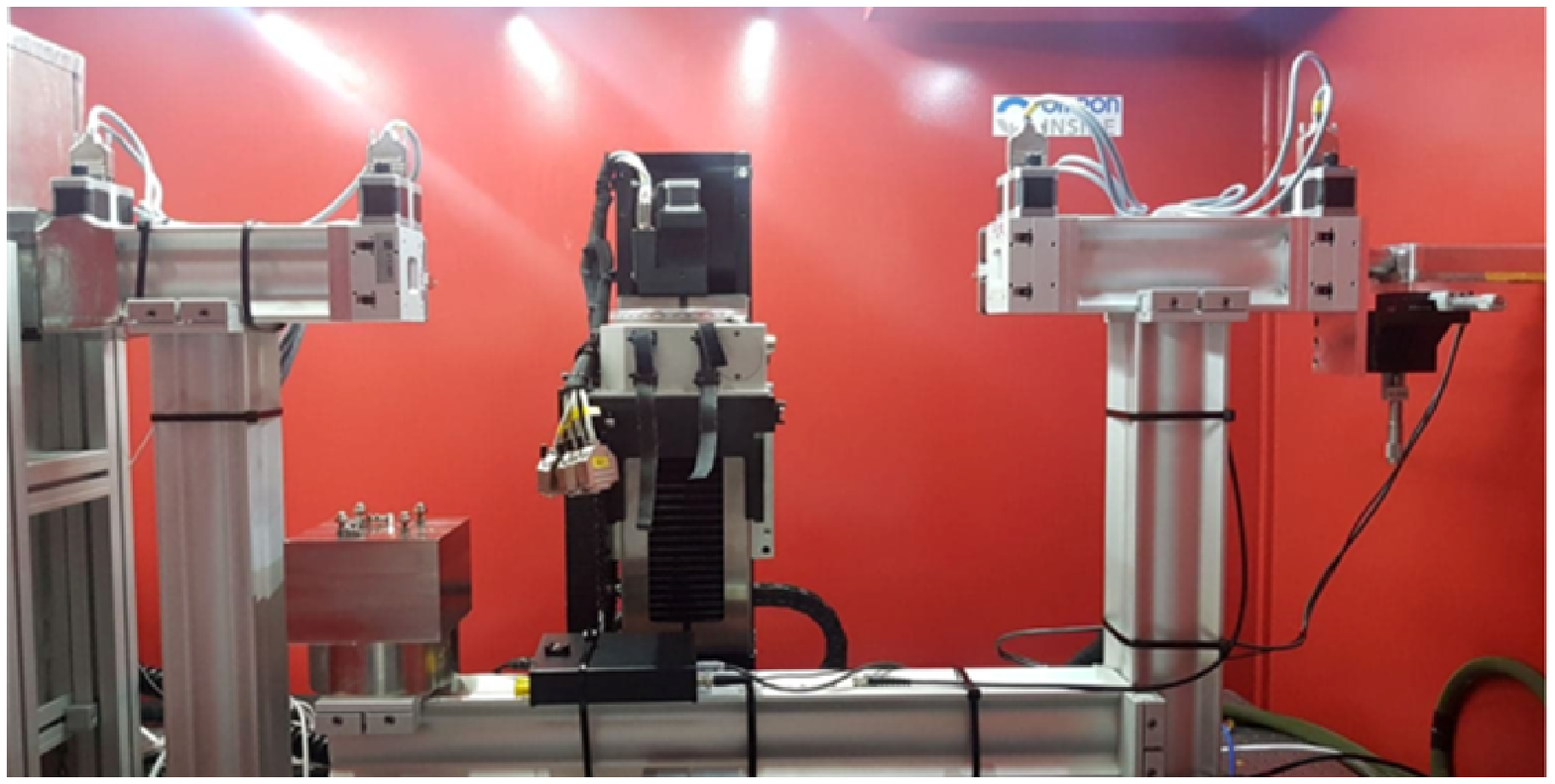}}
	\hspace{.0em}
	{\includegraphics[width=.325\textwidth]{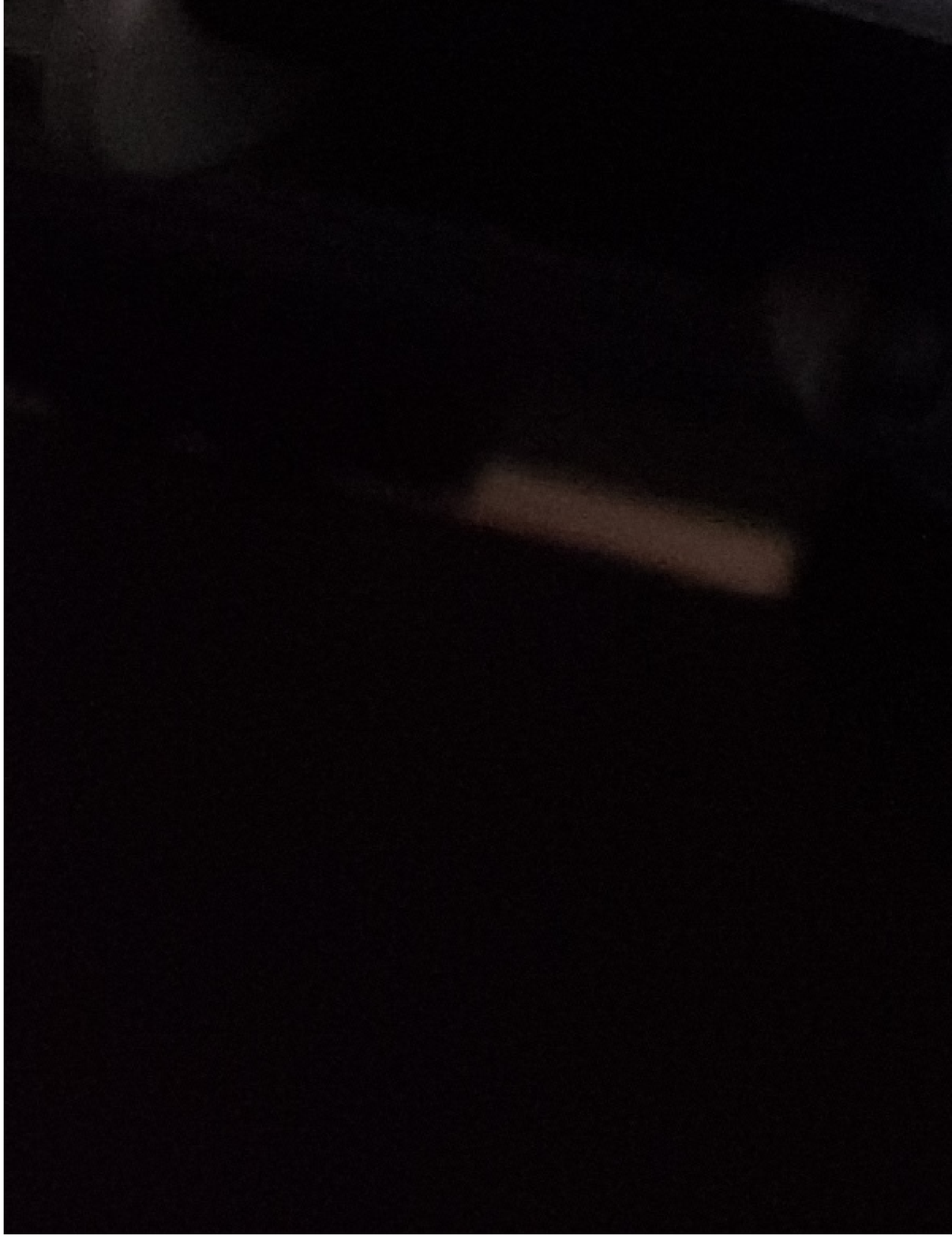}}
	\caption{(a) Scanner setup and (b) support specimen for the cylinders of different materials.}\label{fig:lab}
\end{figure}

We have assessed the algorithm with a real acquisition using a photon counting detector CT scanner depicted in Figure \ref{fig:lab}(a) with a specimen in (b) containing cylinders of iodine and PMMA at different concentrations. The geometry is cone beam with source to axis of rotation distance of 8.1 cm and source to origin distance of 11.4 cm.

\begin{figure}[!h]
	\centering
	{\includegraphics[width=.37\textwidth]{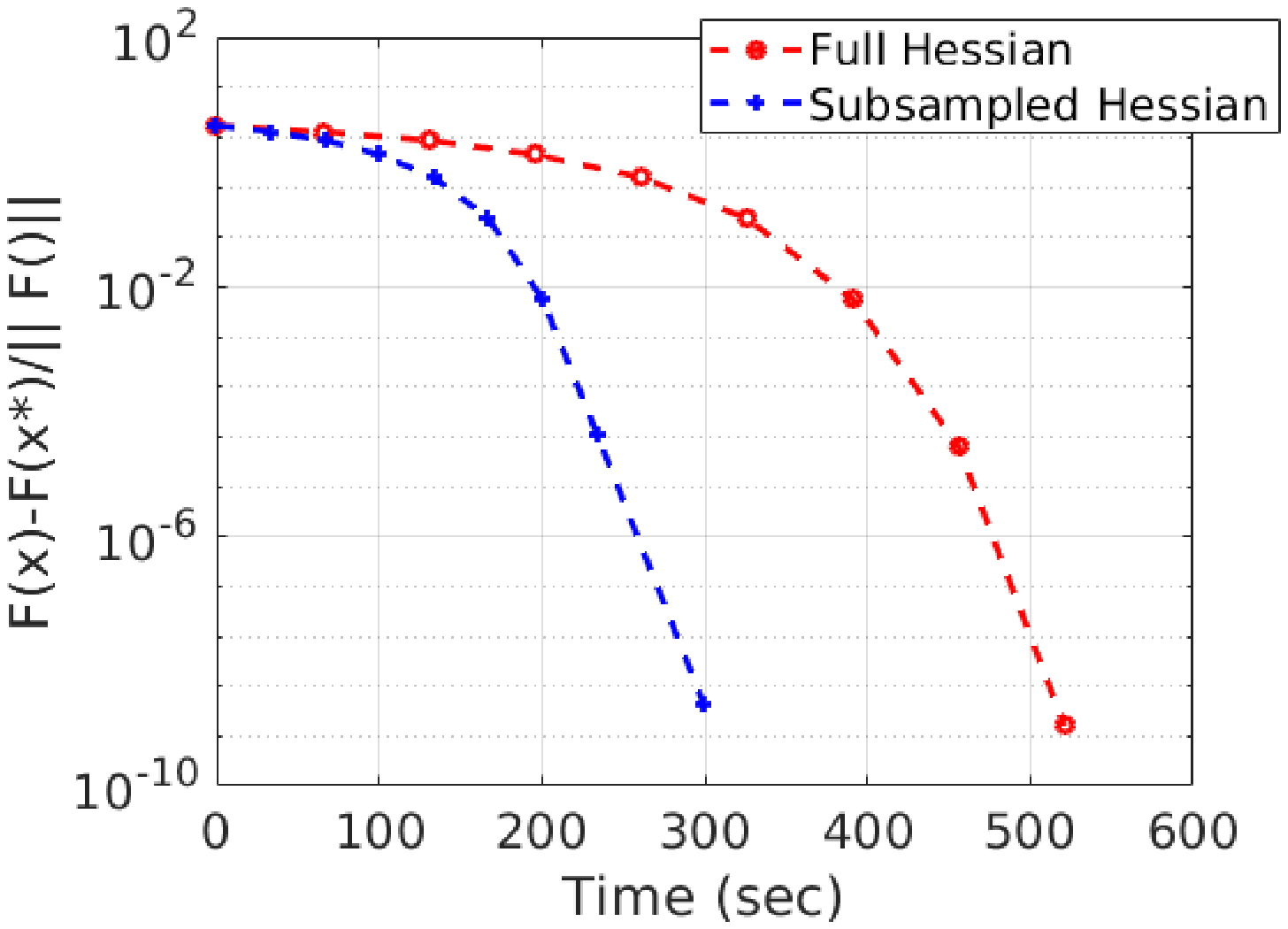}}
	\hspace{.0em}
	{\includegraphics[width=.27\textwidth]{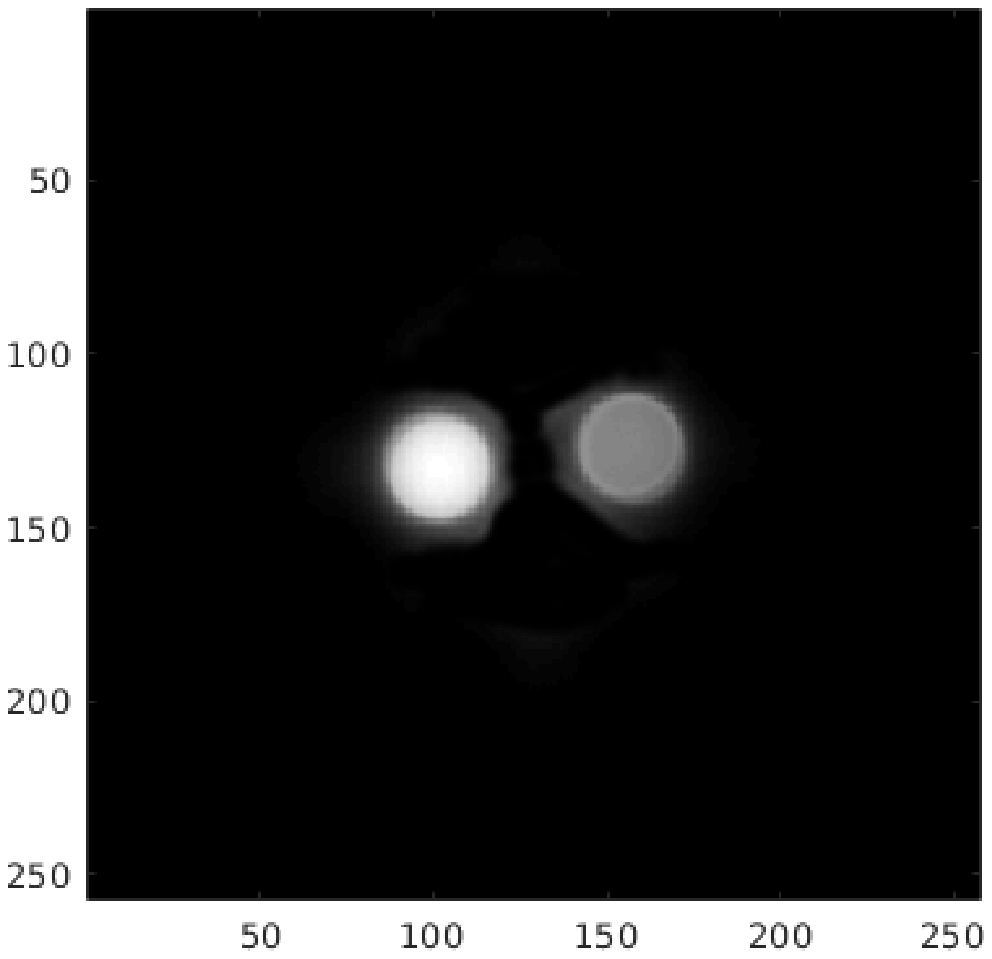}}
	\hspace{.0em}
	{\includegraphics[width=.27\textwidth]{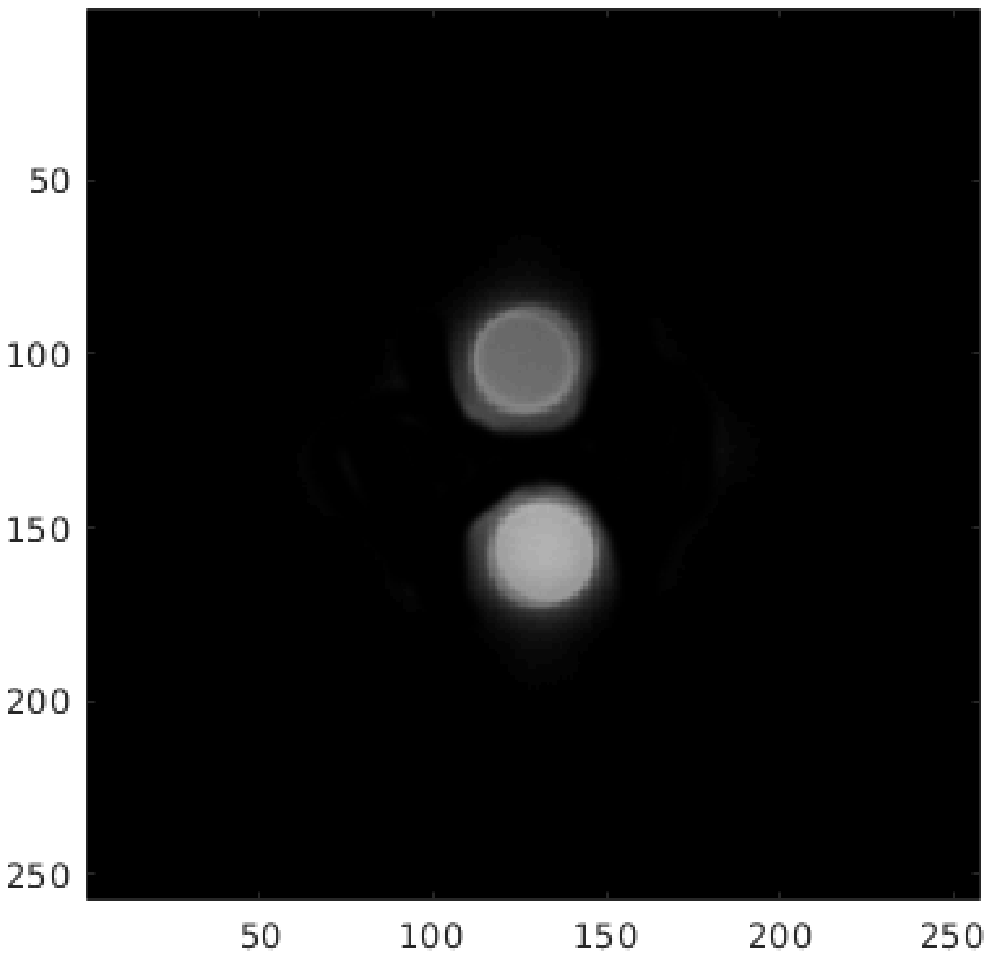}}
	\caption{Error in objective function and reconstruction of the cylinders phantoms.}\label{fig:labp}
\end{figure}

Figure \ref{fig:labp}(a) show the decrease of normalized cost function for the Denoising-IHS subsamplig case and full Hessian; it is worth noting the reduction in computation although is slower compared to the results achieved in the numerical simulations. Table (\ref{tab:error}) and Figures \ref{fig:labp}(b-c) show the quantitative and qualitative error in concentration of the iodine.    

\begin{table}[!h]
	\caption{Iodine concentration error}\label{tab:error}
	\centering
	\begin{tabular}{c|c}
		Ground truth [mg/ml] & estimated \\\hline
		8 & 8.15 \\
		16 & 15.72 \\
	\end{tabular}
\end{table}

\section{Conclusion}
We propose a preliminary study for an efficient second order optimization method, Denoising-IHS, for multi-material decomposition in spectral CT which combines dimensionality reduction through Hessian matrix non-uniform random sub-sampling and regularization by denoising. A practical procedure for efficiently calculating the non-uniform sampling probabilities is shown which involves Fourier spectral approximation for estimating the ridge leverage scores. Furthermore, we analyze the CG solver which does not require the computation of the full Jacobian of the denoiser but only directional derivatives with a substantial reduction in computation. This problem is tested on numerical spectral CT simulations with Poisson noise model. We evaluate the convergence of the algorithm and we are able to show improvements both in terms of qualitative and quantitative accuracy compared to state-of-the-art algorithms. Furthermore, we show reduction in the computational complexity respect to the Newton method involving the full Hessian matrix of the loss function. Further analysis to assess the proposed framework on more complex phantom with experimental acquisitions will be conducted together with a tighter analysis on the lower bounds for the dimensional reduction that it is possible to theoretically achieve and an analysis of the robustness of the algorithm to training sample dimension.

\vskip6pt

\enlargethispage{20pt}

\section*{Acknowledgments}

\noindent The research leading to these results has received funding from the European Union's Horizon 2020 research and innovation programme under the Marie Sklodowska-Curie grant agreement no. 713683 (COFUNDfellowsDTU).
\noindent We would like to sincerely thank Jan Kehres for providing the experimental spectral CT dataset. 
The experimental spectral CT dataset used in Section \ref{sec:exp_res} is available in \cite{Perelli2021data}.

\appendix

\section{Convergence Analysis}\label{app:a2}
	
Recalling the definition of the function $g(\*x) = f(\*x) + \rho_{\nu}(\*x)$ that we aim at minimizing as in (\ref{eq:pr_GLM}) where we hide for simplicity the dependency over $\*y$, we use the denotation  $\Bar{\*H}(\*x)=\nabla^2 g(\*x)$ to define the full Hessian, i.e. the one that would be obtained without block sub-sampling, and $\*H^t(\*x^t)$ the block sub-sampled Hessian defined in (\ref{eq:NewtonSketch}). We assume the following conditions  
\begin{assumption}
The Hessian of the cost function $g(\*x)$ is Lipschitz continuous, therefore there is there exists a constant $L\geq 0$ such that
\begin{equation}\label{eq:Lcont}
\|\Bar{\*H}(\*w) - \Bar{\*H}(\*z) \| \leq L\| \*w - \*z \|, \quad \forall \*w, \*z \in \R^n
\end{equation}
\end{assumption}
	
\begin{assumption} 
For any set $S_t$ with $|S_t|=s$, there is a positive constant such that
\begin{equation}
    \lambda_{min}(\*H^t(\*x^t)) \*I \;\preceq\; \*H^t(\*x^t) \;\preceq\; \lambda_{max}(\*H^t(\*x^t))\*I
\end{equation}
Assuming $\Bar{\*H}(\*x^t)$ the Hessian that would be obtained without block sub-sampling, the error compared to the Hessian obtained by data-loss term sketching is bounded according to the following condition
\begin{equation}\label{eq:cond}
\| \*H^t(\*x^t) - \Bar{\*H}(\*x^t) \| \leq \epsilon \| \Bar{\*H}(\*x^t)   \| 
\end{equation}
	
\noindent Given (\ref{eq:cond}) we have 
\begin{eqnarray}
\lambda_{min}(\*H^t(\*x^t)) & \geq & \lambda_{min}(\Bar{\*H}(\*x^t)) - \epsilon \lambda_{max}(\Bar{\*H}(\*x^t)) \;\; = \;\; \left(1 - \epsilon \kappa (\Bar{\*H}(\*x^t)) \right) \lambda_{min}(\Bar{\*H}(\*x^t)) \nonumber \\
& \geq & \left(1 - 2\epsilon \kappa \right) \lambda_{min}(\Bar{\*H}(\*x^t)) = \lambda_{b}(\Bar{\*H}(\*x^t))
\end{eqnarray}
where $\kappa$ is the condition number of the problem. 
\end{assumption}

We consider the expected error between the oracle solution and the estimate at iteration $t$ where the expectation is taken over the random block selection at iteration $t$
\begin{eqnarray}\label{eq:err_converg}
\Expect_t\left\| \*x^{t+1} - \*x_* \right\| & = & \Expect_t\left\| \*x^t  + \*p_r^t - \*x_* \right\|  \\
& \leq & \underbrace{\Expect_t\left\| \*x^t - 
\left[\*H^t(\*x^t)\right]^{-1} \nabla g(\*x^t) - \*x_*  \right\|}_{\text{$1^{st}$ term}} + 
\underbrace{\Expect_t\left\|  \*p_r^t + \left[\*H^t(\*x^t)\right]^{-1} \nabla g(\*x^t) \right\|}_{\text{$2^{nd}$ term}} \nonumber
\end{eqnarray}
	
\noindent Regarding the $1^{st}$ term in (\ref{eq:err_converg}), we have that 
	
\begin{eqnarray*}
& & \Expect_t\left\| \*x^t - \left[\*H^t(\*x^t)\right]^{-1} \nabla g(\*x^t) - \*x_*  \right\|  =  \Expect_t\left\| \left[\*H^t(\*x^t)\right]^{-1} \left( \*H^t(\*x^t) (\*x^t - \*x_*)  - \nabla g(\*x^t) \right) \right\| \nonumber \\
& \leq & \frac{1}{\lambda_b(\Bar{\*H}(\*x^t))} \Expect_t\left\| \left( \*H^t(\*x^t) - \Bar{\*H}(\*x^t) \right) (\*x^t - \*x_*) + \Bar{\*H}(\*x^t) (\*x^t - \*x_* ) - \nabla g(\*x^t)  \right\| \nonumber \\ 
& \leq & \frac{1}{\lambda_b(\Bar{\*H}(\*x^t))} \left[\underbrace{\Expect_t\left\| \left( \*H^t(\*x^t) - \Bar{\*H}(\*x^t) \right) (\*x^t - \*x_*) \right\|}_{\text{term (a)}} + \underbrace{ \left\| \nabla^2 g(\*x^t) (\*x^t - \*x_* ) - \nabla g(\*x^t) \right\| }_{\text{term (b)}} \right]
\end{eqnarray*}
	
\noindent The term (a) represents the error due to the Hessian sketching and using condition (\ref{eq:cond}) we have
\begin{equation*}
\Expect_t\left\| \left( \*H^t(\*x^t) - \Bar{\*H}(\*x^t) \right) (\*x^t - \*x_*) \right\| \;\; \leq \;\; \epsilon \|\Bar{\*H}(\*x^t) \| \| \*x^t - \*x_* \| \;\; = \;\; \epsilon \lambda_{max}\left(\Bar{\*H}(\*x^t)\right)\| \*x^t - \*x_* \|  
\end{equation*}

\noindent For the term (b), using the Lipschitz continuity of the Hessian (\ref{eq:Lcont}) we have
	
\begin{eqnarray}
\left\| \nabla^2 g(\*x^t) (\*x^t - \*x_* ) - \nabla g(\*x^t) \right\| & \leq & \| \*x^t - \*x_*  \|  \int_{t=0}^1 \| \nabla^2 g(\*x^t) - \nabla^2 g(\*x^t + t(\*x_* - \*x^t))dt \| \nonumber \\
& \leq & \| \*x^t - \*x_*  \|^2  \int_{t=0}^1 L tdt = \frac{L}{2} \| \*x^t - \*x_*  \|^2
\end{eqnarray}
	
\noindent Therefore, the overall $1^{st}$ term is bounded by
\begin{eqnarray}\label{eq:1st_term}
\Expect_t\left\| \*x^t - \left[\*H^t(\*x^t)\right]^{-1} \nabla g(\*x^t) - \*x_*  \right\| & \leq & \frac{1}{(1 - 2\epsilon \kappa) \lambda_{min}(\Bar{\*H}(\*x^t))} \left(\frac{L}{2} \| \*x^t - \*x_* \|^2  + \right. \nonumber \\
& & + \; \epsilon \lambda_{max}(\Bar{\*H}(\*x^t))\| \*x^t - \*x_*  \| \bigg)
\end{eqnarray}
	
\noindent The $2^{nd}$ term in (\ref{eq:err_converg}) represents the residual error after $r$ CG iterations. The convergence rate of the CG method varies at every iteration depending on the spectrum $\{\lambda_1\leq \lambda_2 \leq ,\ldots, \leq \lambda_n\}$ of the positive definite matrix $\*H^t(\*x^t)$. After $r$ steps of the CG method applied to the linear system in (\ref{eq:appr_New}), the iterate $p^t_r$ satisfies
\begin{equation}\label{eq:CG_rate2}
\big\|\*p^t_r - \*p^t_* \big\|^2_{\*C} \leq 
2\sqrt{\kappa(\*C)} \left( \frac{\sqrt{\kappa(\*C)} - 1}{\sqrt{\kappa(\*C)} + 1} \right)^r \big\|\*p^t_0 - \*p^t_* \big\|^2_{\*C}, \quad \*C = \*H^t(\*x^t)
\end{equation}
	
\noindent where $\*p^t_*$ indicates the exact solution and $\| \*x \|^2_{\*C} \triangleq \*x^T\*C\*x$. Using the CG convergence rate (\ref{eq:CG_rate2}) and assuming $\*p^t_0=0$ and $\*p^t_* = - \left[\Bar{\*H}(\*x^t)\right]^{-1} \nabla g(\*x^t)$, we obtain 
\begin{equation}\label{eq:CG_rate_conv}
\left\|  \*p_r^t + \left[\*H^t(\*x^t)\right]^{-1} \nabla g(\*x^t) \right\|_{\*H^t(\*x^t)} \leq 
2\sqrt{\kappa(\*C)} \left( \frac{\sqrt{\kappa(\*C)} - 1}{\sqrt{\kappa(\*C)} + 1} \right)^r \left\| \left[\*H^t(\*x^t)\right]^{-1} \nabla g(\*x^t) \right\|_{\*H^t(\*x^t)} 
\end{equation}
	
\noindent Given $\|\*a \|^2_{\*C} \leq \|\*b \|^2_{\*C}$, we can obtain the following relationship in terms of unweighted norm
\begin{equation}
\lambda_1\|\*a \|^2  \leq \*a^T\*C\*a \leq \*b^T\*C\*b \leq \lambda_n\|\*b \|^2 \Rightarrow \|\*a \| \leq \sqrt{\kappa(\*C)}\|\*b \| \leq \sqrt{\kappa} \|\*b \|
\end{equation}
	
	
	

\noindent By substituting the (\ref{eq:CG_rate_conv}) in the $2^{nd}$ term of (\ref{eq:err_converg}), it results
\begin{eqnarray}\label{eq:CG_res_error}
\Expect_t \left\|\*p_r^t + \left[\*H^t(\*x^t)\right]^{-1} \nabla g(\*x^t) \right\| 
& \leq & 2\sqrt{\kappa(\*C)} \left( \frac{\sqrt{\kappa(\*C)} - 1}{\sqrt{\kappa(\*C)} + 1} \right)^r \Expect_t \left\| \left[\*H^t(\*x^t)\right]^{-1} \nabla g(\*x^t) \right\| 
\nonumber \\
& \leq &  2\sqrt{\kappa(\*C)} \left( \frac{\sqrt{\kappa(\*C)} - 1}{\sqrt{\kappa(\*C)} + 1} \right)^r  \nabla g(\*x^t)\Expect_t \left\| \left[\*H^t(\*x^t)\right]^{-1} \right\| 
\nonumber \\
& \leq &  \frac{2L \sqrt{\kappa(\*C)}}{\mu} \left( \frac{\sqrt{\kappa(\*C)} - 1}{\sqrt{\kappa(\*C)} + 1} \right)^r  \left\| \*x^t -\*x_* \right\|
\end{eqnarray}
	
	

\noindent By substituting the expressions (\ref{eq:1st_term}) and (\ref{eq:CG_res_error}) in (\ref{eq:err_converg}), we have
\begin{eqnarray}\label{eq:final_proof}
\Expect_t\left\| \*x^{t+1} - \*x_* \right\| 
& \leq & \frac{1}{\left(1 - 2\epsilon \kappa \right) \lambda_{min}(\Bar{\*H}(\*x^t))} \left(\frac{L}{2} \| \*x^t - \*x_* \|^2 + \right.\nonumber \\
& & + \;\epsilon \lambda_{max}(\Bar{\*H}(\*x^t))\| \*x^t - \*x_*  \| \bigg) + 2\sqrt{\kappa(\*C)} \left( \frac{\sqrt{\kappa(\*C)} - 1}{\sqrt{\kappa(\*C)} + 1} \right) \left\| \*x^{t} - \*x_* \right\| \nonumber \\
& \leq & \frac{L}{2\left(1 - 2\epsilon \kappa \right) \lambda_{min}(\Bar{\*H}(\*x^t))} \left\| \*x^t - \*x_* \right\|^2 + \nonumber \\
& & + \;\frac{\epsilon \lambda_{max}(\Bar{\*H}(\*x^t))}{2\left(1 - 2\epsilon \kappa \right) \lambda_{min}(\Bar{\*H}(\*x^t))}  \left\| \*x^t - \*x_*  \right\| + 2\sqrt{\kappa(\*C)} \left( \frac{\sqrt{\kappa(\*C)} - 1}{\sqrt{\kappa(\*C)} + 1} \right) \left\| \*x^{t} - \*x_* \right\| \nonumber
\end{eqnarray}



\bibliographystyle{IEEEtran}
\bibliography{biblio}

\end{document}